\pdfoutput=1
\RequirePackage{ifpdf}
\ifpdf 
\documentclass[pdftex]{sigma}
\else
\documentclass{sigma}
\fi

\usepackage{tikz}
\usetikzlibrary{arrows}

\numberwithin{equation}{section}

\newtheorem{Theorem}{Theorem}[section]
\newtheorem*{Theorem*}{Theorem}
\newtheorem{Corollary}[Theorem]{Corollary}
\newtheorem{Lemma}[Theorem]{Lemma}
\newtheorem{Proposition}[Theorem]{Proposition}
\newtheorem{Conjecture}[Theorem]{Conjecture}

\theoremstyle{definition}

\newtheorem{Example}[Theorem]{Example}
\newtheorem{Remark}[Theorem]{Remark}

\begin{document}

\newcommand{\arXivNumber}{2203.15642}

\renewcommand{\thefootnote}{}

\renewcommand{\PaperNumber}{087}

\FirstPageHeading

\ShortArticleName{Generalized $q$-MZVs and Characters of Vertex Algebras}

\ArticleName{Generalized $\boldsymbol{q}$-MZVs and Characters\\ of Vertex Algebras\footnote{This paper is a~contribution to the Special Issue on Recent Advances in Vertex Operator Algebras in honor of James Lepowsky. The~full collection is available at \href{https://sigma-journal.com/Lepowsky.html}{https://sigma-journal.com/Lepowsky.html}}}

\Author{Antun MILAS}

\AuthorNameForHeading{A.~Milas}

\Address{Department of Mathematics and Statistics, University at Albany, State University of New York,\\
1400 Washington Avenue, Albany NY 12222, USA}
\Email{\mail{amilas@albany.edu}}

\ArticleDates{Received December 30, 2025, in final form August 20, 2026; Published online September 02, 2026}

\Abstract{We analyze three related families of $q$-series arising from characters of vertex (super)algebras. More precisely, we consider: (i)~graph series and characters of principal subspaces associated with arc algebras; (ii)~characters, or indices, of class-$\mathcal{S}$ vertex operator algebras in the context of 4d $\mathcal{N}=2$ SCFT; and (iii)~supercharacters of the $\mathcal{U}$ family of vertex superalgebras and their quasi-modular properties. Along the way, we introduce a family of multiple $q$-zeta values associated with simple Lie algebras and present their conjectural
properties.}


\Classification{11F11; 17B69; 11M32; 13F55}

\rightline{\emph{For Jim, on his 80th birthday}}

\renewcommand{\thefootnote}{\arabic{footnote}}
\setcounter{footnote}{0}

\section{Introduction and summary of results}

In this paper we investigate several families of $q$-series arising from characters of vertex algebras and vertex superalgebras, together with Hilbert series of arc algebras. We begin with graph series associated with certain principal subspaces of non-conformal vertex algebras, where modularity is absent but rich combinatorial structures emerge. We then turn to conformal vertex (super)algebras, whose characters are expected or are shown to possess quasi-modular properties. A common theme throughout the paper is the appearance of $q$-multiple zeta values and their generalizations.

There is a vast literature on multiple zeta values (MZVs), much of it devoted to understanding the relations among them. Under the usual multiplication, MZVs (and $q$-MZVs) form an algebra. This property led to the development of the harmonic (or ``stuffle'') and shuffle algebras \cite{Hoffman}. There remain many deep questions concerning generators of these and related algebras \cite{BK2,BK1,Bradley,Goncharov,Zagier1,Ohno,Zagier2,Zudilin2,Zudilin1}, and multiple zeta values have remarkable connections with many areas of mathematics.

In their standard form,\footnote{Another ``standard'' form is $(1-q)^{a_1+ \cdots + a_k}\zeta_q(a_1,\dots,a_k)$ \cite{Bradley}.} $q$-MZVs are usually defined as
\begin{equation} \label{mzvq}
\zeta_q(a_1,\dots,a_k):=\sum_{n_1 > n_2 > \cdots > n_k \geq 1} \frac{q^{(a_1-1)n_1+\cdots + (a_k-1)n_k}}{(1-q^{n_1})^{a_1} \cdots (1-q^{n_k})^{a_k}},
\end{equation}
where $a_i \in \mathbb{N}$ and $a_1 \geq 2$. Throughout the paper, $\mathbb{N}=\{1,2,3,\dots\}$ and $\mathbb{Z}_{\geq 0}=\{0,1,2,\dots\}$.
There are closely related $q$-MZV versions with non-strict summation variables (so-called {\em multiple q-zeta star values}); see \cite{Ohno} and especially \cite{BK2,BK3} for a discussion of various models.

In this work, we are concerned with $q$-MZVs and related $q$-series in connection with characters of vertex (super)algebras.
Throughout the paper we consider vertex algebras, vertex operator algebras,
and vertex operator superalgebras. A vertex operator algebra is a vertex algebra
equipped with a conformal vector, whose associated grading determines the usual notion of character
and is the basis for studying its modular properties. By contrast, a
general vertex algebra need not possess a conformal vector, although it may
still carry natural gradings. In particular, the principal subspaces studied below are naturally graded but are not conformal, so their characters
are defined with respect to this auxiliary grading. We also consider vertex (operator) superalgebras, the
$\mathbb Z_2$-graded analogues of vertex (operator) algebras; whenever the distinction is
irrelevant, we use the collective term ``vertex (super)algebra''.

Our
approach is motivated on the one hand by our previous work on arc spaces and graph series \cite{BJM,JM1,JM2,HM} (see also \cite{Hao}), and on the other hand by character formulas for certain vertex algebras introduced by Arakawa \cite{Arakawa}, related to genus-zero class-$\mathcal{S}$ theories in physics \cite{BR,Beem}.
As already pointed out in \cite{BJM,JM2,Hao,HM}, the Hilbert series of the arc algebra $J_\infty(R_\Gamma)$, where $R_\Gamma$ is the edge algebra (defined below) of a simple\footnote{It is possible to slightly relax this condition and consider graphs with loops~\cite{Hao}.} graph $\Gamma$, coincides with the character of a certain {\em principal subspace} $W_\Gamma$ (after Feigin and Stoyanovsky~\cite{FS}) inside an appropriate lattice vertex algebra \cite{CLM,Kawasetsu, MP}. Moreover, the arc algebra $J_\infty(R_\Gamma)$ carries the structure of a commutative vertex algebra, and this correspondence can be realized as an isomorphism of commutative graded vertex algebras. We stress, however, that the principal subspace $W_\Gamma$ is not a conformal vertex algebra, so its grading is not canonical.
Therefore, following \cite{Hao}, we equip the principal subspace $W_\Gamma$ with the grading
\[
\deg \bigl({\rm e}^{\alpha_i}_{(-m-1)}\bigr)=m+b_i,
\]
where
\[
Y({\rm e}^{\alpha_i},x)=\sum_{m\in\mathbb{Z}} {\rm e}^{\alpha_i}_{(m)}x^{-m-1},
\]
together with the color grading recording the multiplicity of each generator ${\rm e}^{\alpha_i}$. This grading is the analogue of the grading on the arc algebra $J_\infty(R_\Gamma)$, where the jet variables satisfy
\[
\deg(x_{i,m})=m+b_i.
\]
Thus the correspondence
\[
x_{i,m}
\longleftrightarrow
{\rm e}^{\alpha_i}_{(-m-1)}
\]
preserves both the homogeneous grading and the color grading.

Unless stated otherwise, the variable $q$ records the chosen grading
(conformal in the conformal case and the auxiliary grading described above in
the non-conformal case), while the variables~$x_i$ record the multigrading (or
charge) associated with the generators.

We now return to the connection with arc algebras.
The corresponding Hilbert series of $J_\infty(R_\Gamma)$ is what we call the {\em graph series}
\begin{equation*}
H_{\Gamma,{\bf b}}(q):=\sum_{{\bf n} \in\mathbb{N}_0^r}
\frac{q^{\frac12 {\bf n}A_\Gamma{\bf n}^{\mathsf{T}}+b_1n_1+\cdots+b_rn_r}}
{(q)_{n_1}\cdots(q)_{n_r}},
\end{equation*}
where $A_\Gamma$ is the (symmetric) adjacency matrix of a graph $\Gamma$ with $r$ vertices and ${\bf b}=(b_1,\dots,b_r)\in\mathbb{N}^r$.
Here we use the standard $q$-Pochhammer notation
\[
(q)_n=\prod_{j=1}^n(1-q^j) \qquad \text{and} \qquad (q)_\infty=\prod_{j=1}^\infty(1-q^j),
\]
 whenever the latter product is viewed formally or converges analytically.
The vector ${\bf b}$ determines the grading of the generators of $W_\Gamma$. Equivalently, we can start from the {\em full} character of~$W_\Gamma$,
\begin{equation*}
H_{\Gamma}(q,{\bf x}):=
\sum_{{\bf n}\in\mathbb{N}_0^r}
\frac{{\bf x}^{\bf n}\,
q^{\frac12{\bf n}A_\Gamma{\bf n}^{\mathsf{T}}}}
{(q)_{n_1}\cdots(q)_{n_r}},
\end{equation*}
and specialize to ${\bf x}=q^{\bf b}$. Essentially the same type of $(q,x)$-series appeared in so-called ``knots-quivers'' correspondences; see \cite{Ekholm,Quivers}.

In this paper, we are only interested in the standard grading ${\bf b}=(1,\dots,1)$, in which case
\[
\deg \bigl({\rm e}^{\alpha_i}_{(-m-1)}\bigr)=m+1,
\]
and we simply write $H_\Gamma(q)$ instead of $H_{\Gamma,{\bf b}}(q)$, thus recovering the original graph series introduced in \cite{BJM,JM1,JM2}. Although not directly related to this work, it is worth noting that there is an interesting connection among the Donaldson--Thomas invariants of symmetric quivers, cohomological Hall algebras, and the dual of the principal subspace of $W_{\Gamma}$ \cite{DFR,DM} (see also~\cite{Efimov}).

In parallel with the $h$-polynomial of the Hilbert series of a finitely generated graded commutative algebra, it is convenient to
express $H_{\Gamma}(q)$ in the form
\[
H_{\Gamma}(q)=\frac{h_\Gamma(q)}{(q)_\infty^d},
\]
where $d$ is the Krull dimension of the {\em edge algebra} $R_\Gamma$. Here \begin{equation} \label{gamma}
R_\Gamma:=k[x_1,\dots,x_r]/I_\Gamma,
\end{equation}
 where $I_\Gamma$ is the ideal generated by all monomials $x_i x_j$ such that the vertices $i$ and $j$ of $\Gamma$ are connected by an edge; there are no other generators.

In our previous paper with Bringmann and Jennings-Shaffer \cite{BJM}, see also \cite{Hao,HM}, we had observed that graph series of several well-known graphs (with a small number of vertices) give rise to interesting $q$-series with peculiar modular and combinatorial properties. The problem of finding a closed expression for $H_\Gamma(q)$ for every graph is of course hopeless, but in many examples we have obtained significant simplification that led to new $q$-series identities.

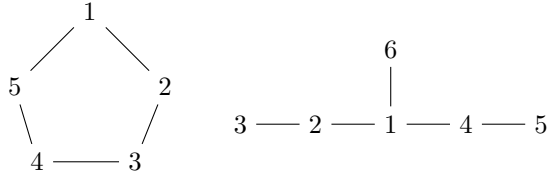
\begin{figure}[!ht]\centering
\begin{tikzpicture}[scale=1]
 \node (n1) at (9,3.5) {1};
 \node (n3) at (7,3.5) {3};
 \node (n2) at (8,3.5) {2};
 \node (n4) at (10,3.5) {4};
 \node (n5) at (11,3.5) {5};
 \node (n6) at (9,4.5) {6};
 \node (n10) at (4.3,3) {4};
 \node (n11) at (4,4) {5};
 \node (n7) at (5,5) {1};
 \node (n8) at (6,4) {2};
 \node (n9) at (5.6,3) {3};
 \foreach \from/\to in {n3/n2,n2/n1,n1/n4,n4/n5,n1/n6,n7/n8,n8/n9,n9/n10,n10/n11,n11/n7}
 \draw (\from) -- (\to);
\end{tikzpicture}

\caption{$C_5$ and $E_6$ graphs.}\label{Figure1}
\end{figure}

We noticed in \cite{BJM} that non-isomorphic graphs can have identical graph series. For example, the $5$-cycle $C_5$ (later also denoted by $\Gamma_5$) and $E_6$
graphs (see Figure~\ref{Figure1})
give the same graph series up to an Euler factor, i.e., a multiplicative power of the Euler product $(q)_\infty$:
\begin{equation} \label{c5e6}
H_{C_5}(q)=q^{-1} \frac{\sum_{n \geq 1} \frac{n q^n}{1-q^n}}{(q)_\infty^2}, \qquad H_{E_6}(q)=q^{-1} \frac{\sum_{n \geq 1} \frac{q^n}{(1-q^n)^2}}{(q)_\infty^3},
\end{equation}
and thus
\begin{equation} \label{c5-e6}
H_{C_5 \oplus {\rm pt} }(q)=H_{E_6}(q),
\end{equation}
where ``${\rm pt}$'' denotes the trivial (or singleton) graph with one vertex and no edges, and~$\oplus$ denotes the disjoint union (coproduct) of graphs.
This indicates that $H_{\Gamma}$, which is clearly a graph invariant, does not distinguish non-isomorphic graphs. This is perhaps not so surprising, as there are non-isomorphic algebraic schemes with identical Hilbert series.

Our first result in this paper extends \eqref{c5e6}--\eqref{c5-e6} to two infinite families
of graphs. Part (i) identifies one family with the special multiple
$q$-zeta values $\mathfrak z_q^*(2,1,\dots,1)$ (see \eqref{ohno} below), while part (ii) recovers the
ordinary $q$-zeta values $\mathfrak z_q(k)$.

\begin{Theorem} \label{main}
We have
\begin{itemize}\itemsep=0pt
\item[$(i)$]
Let $\Gamma_{3k+2}$ be the leafless graph defined in Section~{\rm \ref{section2.1}}. Then
\[
H_{\Gamma_{3k+2}}(q)=q^{-1}\dfrac{\sum_{n \geq 1} \dfrac{n \left(\begin{smallmatrix}n+k-2\\k-1\end{smallmatrix}\right) q^n}{1-q^n}
}{(q)_\infty^{k+1}}.
\]
Here the numerator is the logarithmic derivative of MacMahon's $(k+1)$-dimensional generating series for higher-dimensional partitions {\rm \cite{Andrews0}}:
\[
\frac{1}{\prod_{n \geq 1} (1-q^n)^{{n+k-2 \choose k-1}}}.
\]
\item[$(ii)$]
Let $T_{2k+2}$, $k \geq 2$, be the $(k+1)$-star graph described in Section~{\rm \ref{section2.2}}. Then
\[
H_{T_{2k+2}}(q)=q^{-1}\frac{\mathfrak{z}_{q}(k)}{(q)_\infty^{k+1}},
\]
where $\mathfrak{z}_{q}(k)=\sum_{n \geq 1} \frac{q^n}{(1-q^n)^k}=\sum_{n \geq 1} \frac{{n+k-2 \choose k-1}q^n}{1-q^n}$ is the $k$-th $q$-zeta value.
\end{itemize}
\end{Theorem}
From this, we obtain the following.

\begin{Corollary} For $k \geq 1$,
\[
H_{\Gamma_{3k+2} \oplus pt}(q)=k \cdot H_{T_{2k+4}}(q)-(k-1)H_{T_{2k+2} \oplus pt }(q).
\]
In particular, for $k=1$ this is equivalent to~\eqref{c5-e6}.
\end{Corollary}
In addition to results about graph series, we also analyze Hilbert series of the edge ideals $R_{\Gamma_{3k+2}}$ and of $R_{T_k}$.

Since we have encountered a $q$-zeta–type value $\mathfrak{z}_q(k)$, it is natural to ask whether other $q$-multiple zeta values also arise as graph series, possibly with additional Euler factors.
As indicated above, we are primarily interested in the standard grading, corresponding to $b_i = 1$.
Within this framework, the $q$-MZVs relevant to our work were introduced by Ohno, Okuda, and Zudilin in \cite{Ohno}:
\begin{equation} \label{ohno}
\mathfrak{z}_q^*(a_1,\dots,a_k):=\sum_{n_1 \geq n_2 \geq \cdots \geq n_k \geq 1} \frac{q^{n_1}}{(1-q^{n_1})^{a_1} \cdots (1-q^{n_k})^{a_k}}.
\end{equation}
The star symbol indicates that the summation is over non-strict summation variables as opposed to the more standard model:
\[
\mathfrak{z}_q(a_1,\dots,a_k):=\sum_{n_1 > n_2 > \cdots > n_k \geq 1} \frac{q^{n_1}}{(1-q^{n_1})^{a_1} \cdots (1-q^{n_k})^{a_k}}.
\]
The difference between the strict and non-strict $q$-MZV models is controlled by the ``lower depth'' $q$-MZVs, so one can easily transform between the two.
Our next result concerns graph series and multiple $q$-zeta values (\ref{ohno}).
\begin{Theorem} \label{qmzv} For every choice of positive integers $a_1,\dots,a_k$, there is a simple graph $Z_{a_1,\dots,a_k}$ such that
\[
H_{Z_{a_1,\dots,a_k}}(q)=\frac{q^{-1} \mathfrak{z}_q^*(a_1,\dots,a_k)}{(q)^{k+a_1+\cdots +a_k}_\infty}.
\]
In particular, for $k=1$ and $a_1=r$, this recovers $H_{T_{2r+2}}(q)$ and Theorem~{\rm \ref{main}}\,$(ii)$.
\end{Theorem}

\begin{Remark}
{ With a different framing vector ${\bf b} \in \mathbb{N}$ and the same underlying graph as in Theorem~\ref{qmzv}, the graph series
of $Z_{a_1,\dots,a_k}$ can be expressed using the standard $q$-MZVs (\ref{mzvq}), albeit in non-strict form.}
\end{Remark}

Our previous results concerned non-conformal vertex algebras, and therefore their modular properties cannot be studied using standard methods such as those of~\cite{Zhu}.
In Section~\ref{section3}, we change perspective and consider characters of vertex operator algebras arising from the genus zero class-$\mathcal{S}$ vertex operator algebras~\cite{Arakawa,Beem}, which play an important role in four-dimensional $\mathcal{N}=2$ superconformal field theories~\cite{BR,Beem}.
These vertex operator algebras are expected to be \emph{quasi-lisse}, and consequently their characters should satisfy modular linear differential equations (MLDEs) with modular coefficients; see \cite{AK} for a precise formulation.

As we shall see later, certain characters of class-$\mathcal{S}$ vertex operator algebras can be expressed as suitable modifications of $q$-multiple zeta values associated with a simple Lie algebra, which we introduce next.

Denote by $\Delta$ a root system of ADE type, by $\Delta_+$ the set of positive roots, and let $\langle \cdot, \cdot \rangle$ denote the inner product normalized
such that $\langle \alpha,\alpha \rangle=2$ for every root $\alpha$. For integers $k_\alpha \geq 1$, we define a $q$-MZV associated with $\mathfrak{g}$ as
\begin{equation} \label{qmzvg}
\zeta_{\mathfrak{g},q}(k_\alpha; \alpha \in \Delta_+):=\sum_{\lambda \in P_+} \frac{q^{\frac12 \sum_{\alpha \in \Delta_+} k_\alpha \langle\lambda+\rho,\alpha \rangle}}{\prod_{\alpha \in \Delta_+} (1-q^{\langle \lambda, \alpha+\rho\rangle})^{k_\alpha}},
\end{equation}
where the summation is over the cone of positive dominant integral weights. The asymptotic behavior motivates our definition. Indeed, after setting $q={\rm e}^{2\pi {\rm i}\tau}$ and taking the asymptotic limit $\tau \to 0+$ with
$k_{\alpha}:=k \geq 2$, we recover
Witten's $\zeta$-function of $\mathfrak{g}$ (see \cite{Zagier})
\[
\sum_{\lambda \in P_+} \frac{1}{\dim (V_\lambda)^k},
\]
 as the leading coefficient of the asymptotic expansion. This fact follows from the Weyl dimension formula.

In the simplest case of $\mathfrak{g}=\mathfrak{sl}_2$, where $\Delta=\{ \pm \alpha \}$, these $q$-zeta values are given by
\[
\zeta_{\mathfrak{sl}_2,q}(k)=\sum_{n \geq 1} \frac{q^{\frac{k}{2}n}}{(1-q^n)^k}.
\]
This slightly non-standard $q$-zeta value has several nice properties. For instance, the summand is invariant under $q \to q^{-1}$ for $k$ even and skew-invariant for $k$ odd. For simplicity, we only consider $k$ even here.
Then it is known (or it can be easily verified) that for $k \geq 1$, $\zeta_{\mathfrak{sl}_2,q}(2k)$, up to an additive constant, is a linear combination of Eisenstein series, and in fact we have
\begin{equation*} 
\mathbb{C}[\zeta_{\mathfrak{sl}_2,q}(2k): k \geq 1]=\underbrace{\mathbb{C}[E_2(q),E_4(q),E_6(q)]}_{:=\mathcal{QM}},
\end{equation*}
where $E_{2i}$ are Eisenstein series defined with their $q$-expansion; see formula (\ref{eis}) in Section~\ref{section3}.
The space $\mathcal{QM}$ is called the ring of (holomorphic) quasi-modular forms on ${\rm SL}(2,\mathbb{Z})$.
For $k$ odd and $k \geq 3$, these values
capture Eisenstein series of level two which we do not consider in this paper.
In parallel with other, more familiar examples of $q$-MZVs, we expect a symmetry formula to hold.
\begin{Conjecture} \label{mzvgq} Assume that $\mathfrak{g}$ is of~$A$ type. Denote by $\Delta_+=\{ \alpha_1,\dots,\alpha_r \}$ the set of positive roots of~$\mathfrak{g}$.
Then we have
\begin{enumerate}\itemsep=0pt
\item[$(i)$] Let $k_\alpha \in \mathbb{N}$ for every $\alpha \in \Delta_+$. Then
\[
\sum_{\sigma \in W} \zeta_{\mathfrak{g},q}(2k_{\sigma \alpha_1},\dots,2k_{\sigma \alpha_r}) \in \mathcal{QM},
\]
where $W$ is the Weyl group of $\mathfrak{g}$ and where we let $k_{-\alpha_i}=k_{\alpha_i}$.

\item[$(ii)$] In particular, for $k=k_\alpha$, we have
\begin{equation*} 
\zeta_{\mathfrak{g},q}(2k):=\zeta_{\mathfrak{g},q}(2k,2k,\dots,2k) \in \mathcal{QM}.
\end{equation*}
\end{enumerate}
\end{Conjecture}

The second part of the conjecture is closely related to another quasi-modularity conjecture for the characters of vertex operator algebras ${\bf V}_{\mathfrak{g},k}$ introduced by Arakawa \cite{Arakawa}.
Based on previous work of Beem and Rastelli \cite{Beem}, dealing with the $\mathfrak{g}=\mathfrak{sl}_2$ case, we expect the following.

\begin{Conjecture} \label{Arakawa} For every even $k >3$,
\[
\eta(\tau)^{m_{\mathfrak{g},k}} {\rm ch}[{\bf V}_{\mathfrak{g},k}](\tau) \in \mathcal{QM},
\]
where $m_{\mathfrak{g},k}$ is a positive integer depending on $\mathfrak{g}$ and $k$ and $\eta(\tau)=q^{1/24}\prod_{i \geq 1}(1-q^i)$ is the Dedekind eta function.
\end{Conjecture}

In Section~\ref{section3}, using results on the constant terms of quasi-Jacobi forms, we prove (see Theorem~\ref{sl3-qm} and Proposition~\ref{sl3-qm2}) the conjectures in the first non-trivial case:
\begin{Corollary} Conjectures~{\rm \ref{mzvgq}} and~{\rm \ref{Arakawa}} hold for $\mathfrak{g}=\mathfrak{sl}_3$.
\end{Corollary}

In our subsequent work~\cite{BVM}, we established significantly more general
quasi-modularity results related to the conjectures above. Although~\cite{BVM} treats a much broader class of examples, the methods developed there
are substantially different from those used in the present paper. The
$\mathfrak{sl}_3$ case treated here therefore remains an independent and
self-contained proof.

The third part of the paper concerns a different but related family of examples. Motivated by a~similar circle of ideas, Adamovi\'c and the author constructed in \cite{AM} a new family of vertex operator (super)algebras $\mathcal{U}^{(m)}$, $m \geq 3$. More precisely, $\mathcal{U}^{(m)}$ is a $\mathbb{Z}_{\geq 0}$-graded vertex operator superalgebra for odd $m$ and $\frac12 \mathbb{Z}_{\geq 0}$-graded vertex operator algebra for even~$m$.
These are expected to share properties with those in the class-$\mathcal{S}$ setting~\cite{BR, Beem} and with Deligne's series~\cite{AK}, including quasi-modularity of characters.
In contrast to Deligne's series -- where the corresponding vertex operator algebras admit a unique ordinary module -- $U^{(m)}$ (for $m \ge 3$) possesses precisely $m$ inequivalent irreducible ordinary modules. Consequently, their characters are expected to transform under an additional level-$m$ structure.
 The algebra $\mathcal{U}^{(m)}$ is an infinite conformal extension of $L_{-1}(\mathfrak{sl}_m)$ and can also be viewed as a coset subalgebra inside $L_1(\mathfrak{psl}_{m|m})$; see \cite[Proposition~7.1]{AM}. This explains its simultaneous relation to affine and superalgebraic examples.
We finish with a result on supercharacters for the $\mathcal{U}$-series in the vein of Conjecture~\ref{Arakawa}.
\begin{Theorem} For $m$ odd, we have
\[
{\rm sch}[\mathcal{U}^{(m)}](\tau)= \frac{\eta(m \tau)^3}{\eta(\tau)^{2m-1}} \mathbb{F}_m(\tau),
\]
where $\mathbb{F}_m(\tau) \in \mathbb{Q}[E_2(\tau),E_4(\tau),E_6(\tau),E_2(m \tau),E_4(m \tau),E_6(m \tau)]$ is a quasimodular form on $\Gamma_0(m)$.
\end{Theorem}

\section[Graph series and multiple q-zeta values]{Graph series and multiple $\boldsymbol{q}$-zeta values}\label{section2}

This section contains the first part of the paper. We recall the elementary
commutative-algebra background for edge ideals, then compute graph series for
the families $\Gamma_{3k+2}$ and $T_{2k+2}$, and finally explain how the same
construction realizes the non-strict $q$-MZVs appearing in the introduction.
The connection with vertex algebras enters through the interpretation of these
graph series as characters of principal subspaces mentioned in the introduction. However, once this
correspondence is established, the proofs in this section are entirely
commutative-algebraic and combinatorial.

\subsection[Hilbert series of R\_Gamma]{Hilbert series of $\boldsymbol{R_{\Gamma}}$}\label{section2.1}

First we recall a few basic facts about square-free monomial ideals and their Hilbert series.

Let $[n]:=\{1,\dots,n\}$. A simplicial complex $\Delta$ on the vertex set $[n]$ determines the square-free monomial ideal
\[
I_{\Delta}:=\bigl\langle x_{i_1}\cdots x_{i_r}\mid \{i_1,\dots,i_r\}\notin\Delta\bigr\rangle
\subseteq S:=k[x_1,\dots,x_n].
\]
This is the unique \emph{Stanley--Reisner ideal} associated with $\Delta$. Conversely, every square-free monomial ideal $I\subseteq S$ determines a unique simplicial complex by declaring $\{i_1,\dots,i_r\}$ to be a face precisely when $x_{i_1}\cdots x_{i_r}\notin I$. In particular, an edge ideal is square-free because it is generated by square-free quadratic monomials $x_ix_j$ corresponding to the edges of a simple graph.

The Stanley--Reisner ring $k[\Delta]:=S/I_{\Delta}$ has a multigraded Hilbert series, with $\deg (x_i)=1$, that can be computed directly from $\Delta$. We regard it as an element of the formal power-series ring $\mathbb{Z}[[x_1,\dots,x_n]]$. It is given by
\[
\mathrm{H}(k[\Delta]; x_1, \dots, x_n) = \sum_{\sigma \in \Delta} \prod_{i \in \sigma} \frac{x_i}{1 - x_i},
\]
where the sum is over all faces $\sigma$ of $\Delta$ and every denominator is expanded at the origin, namely $(1-x_i)^{-1}=\sum_{r\geq0}x_i^r$. Thus the displayed rational expression is always understood through its image in $\mathbb{Z}[[x_1,\dots,x_n]]$. This formula arises from the fact that the monomials not in the Stanley--Reisner ideal $I_\Delta$ correspond precisely to the monomials supported on faces of $\Delta$, and for each face $\sigma$, the set of monomials supported on $\sigma$ contributes the term $\prod_{i \in \sigma} \frac{x_i}{1 - x_i}$. To obtain the standard (total degree) Hilbert series, we make the substitution $x_i = t$ for all $i$ as our generators all have degree $1$. We denote the corresponding Hilbert series by $H_t(k[\Delta])$.

The combinatorial structure of the simplicial complex $\Delta$ of dimension
$d-1$ is encoded by its $f$-vector
$(f_{-1},f_0,\dots,f_{d-1})$ and its $h$-polynomial
$h(t)=h_0+h_1t+\cdots+h_dt^d$. These are related to the Hilbert
series of the Stanley--Reisner ring $k[\Delta]$ by
\[
H_t(k[\Delta])
=
\sum_{i=0}^{d}\frac{f_{i-1}t^i}{(1-t)^i}
=
\frac{h_0+h_1t+\cdots+h_dt^d}{(1-t)^d}.
\]
In particular, the order of the pole at $t=1$ equals the Krull dimension
$d=\dim k[\Delta]$.

In this work, we are only interested in edge ideals (\ref{gamma}) and their graph series. Since the edge ideal $I_\Gamma=\langle x_ix_j\mid \{i,j\}\in E(\Gamma)\rangle$ of a simple graph is generated by square-free monomials, it is a Stanley--Reisner ideal and $k[\Delta]=R_\Gamma$, where $\Delta$ is the corresponding independence complex. Then $\mathrm{H}(k[\Delta];x_1,\dots,x_n)$ can be
analyzed using the graph series as follows.

\begin{Proposition} Suppose that $\Delta$ is associated with $I_{\Gamma}$; see \eqref{gamma}. Then we have
\[
{H}(k[\Delta];\boldsymbol{x}) ={\rm CT}_q H_{\Gamma}(q; {\bf x}),
\]
where on the right-hand side we are taking the constant term with respect to the $q$-expansion.
\end{Proposition}
\begin{proof} Observe that the contribution to the constant term in $q$, that is monomials of the form $\boldsymbol{x}^{\boldsymbol{n}}=x_1^{n_1} \cdots x_r^{n_r}$, only comes from $r$-tuples ${\bf n}=(n_1,\dots,n_r) \in \mathbb{Z}_{\geq 0}^r $
for which ${\bf n}^{\mathsf{T}} A_{\Gamma} {\bf n}=0$ and that is precisely the case when $x_i x_j$ does not divide $x_1^{n_1} \cdots x_r^{n_r} \notin I$.
\end{proof}

For an edge ideal $I_\Gamma$, the corresponding simplicial complex $\Delta$ and the full Hilbert series $H(R_\Gamma;\boldsymbol{x})$ may be quite complicated to determine explicitly. But in the situation when the
Hilbert series is specialized $x_1=\cdots = x_n=t$, so it only depends on a single variable, we have the following very useful result (see, e.g., \cite{Edge}).

\begin{Lemma} \label{recursion} Let $X_{v}$ be the neighborhood of the node $v$ of the graph $\Gamma$.
Then we have
\[
H_t(R_\Gamma)=H_t(R_{\Gamma \setminus \{v\}})+\frac{t}{1-t}H_t(R_{\Gamma \setminus (\{v\}\cup X_v)}).
\]
\end{Lemma}

Next we prove Theorems~\ref{main} and~\ref{qmzv} stated in the introduction and related results.

\subsection[Edge ideals and graph series of Gamma\_\{3k+2\}, k>=1]{Edge ideals and graph series of $\boldsymbol{\Gamma_{3k+2}}$, $\boldsymbol{k \geq 1}$}\label{section2.2}

We first define a family of leafless graphs $\Gamma_{3k+2}$ that generalize the $5$-cycle and are obtained by gluing together several pentagons. The graph $\Gamma_{3k+2}$ has $3k+2$ vertices, enumerated by $1,2,\dots,3k+2$. The vertices $1$ through $5$ form a pentagon, so we impose the adjacency relations
\[
1 \sim 2 \sim 3 \sim 4 \sim 5 \sim 1.
\]

Continuing, we require
\[
i \sim i+1, \qquad 5 \leq i \leq 3k+1.
\]
Vertices $1$ and $4$ are special: for these we additionally impose
\[
1 \sim 3t+2 \quad \text{and} \quad 4 \sim 3t+2, \qquad t \geq 2.
\]
Finally, we impose the relations
\[
3s+1 \sim 3(\ell+1)+2, \qquad 2 \leq s \leq \ell \leq k-1.
\]

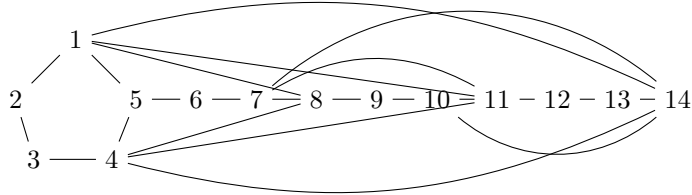
\begin{figure}[!ht]\centering\vspace{-4mm}
\begin{tikzpicture}[scale=0.6]
 \node (n4) at (4.3,2.7) {3};
 \node (n5) at (4,4) {2};
 \node (n1) at (5,5) {1};
 \node (n2) at (6,4) {5};
 \node (n3) at (5.6,2.7) {4};
 \node (n6) at (9,4) {8};
 \node (n7) at (7,4) {6};
 \node (n8) at (8,4) {7};
 \node (n10) at (10,4) {9};
 \node[fill=none] (n9) at (13,4) {11};
\node[fill=none] (n13) at (14.5,4) {12};
 \node (n11) at (11.5,4) {10};
 \node (n14) at (16,4) {13};
 \node (n12) at (17.5,4) {14};
\path (n8) edge [bend left=30] (n9);
\path (n8) edge [bend left=40] (n12);
\path (n1) edge [bend left=20] (n12);
\path (n3) edge [bend right=20] (n12);
\path (n11) edge [bend right=40] (n12);
 \foreach \from/\to in {n4/n5,n5/n1,n1/n2,n2/n3,n3/n4,n2/n7,n7/n8,n8/n6,n3/n6,n1/n6,n6/n10,n10/n11,n11/n9,n9/n13,n1/n9,n3/n9,n13/n14,n14/n12}
 \draw (\from) -- (\to);
\end{tikzpicture}\vspace{-4mm}

\caption{Graph $\Gamma_{14}$.}\label{fig:gamma14}
\end{figure}

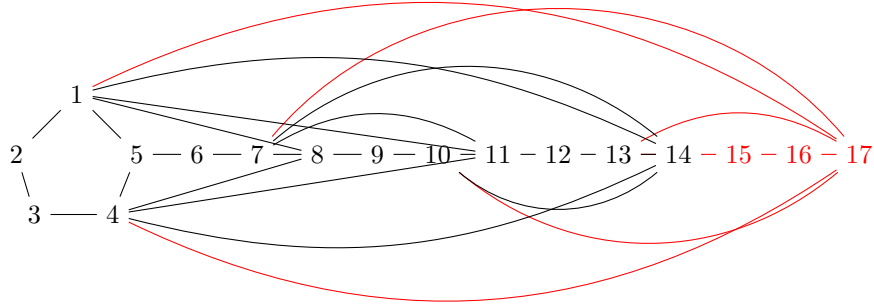
\begin{figure}[!ht]\centering\vspace{-4mm}
\begin{tikzpicture}[scale=0.5]
 \node (n4) at (4.3,2.7) {3};
 \node (n5) at (4,4) {2};
 \node (n1) at (5,5) {1};
 \node (n2) at (6,4) {5};
 \node (n3) at (5.6,2.7) {4};
 \node (n6) at (9,4) {8};
 \node (n7) at (7,4) {6};
 \node (n8) at (8,4) {7};
 \node (n10) at (10,4) {9};
 \node (n11) at (11.2,4) {10};
 \node (n9) at (12.6,4) {11};
 \node (n13) at (14,4) {12};
 \node (n14) at (15.5,4) {13};
 \node (n12) at (17,4) {14};

 \node (n16) at (18.5,4) [color = red] {15};
 \node (n17) at (20,4) [color = red] {16};
 \node (n15) at (21.5,4) [color = red] {17};

\path (n8) edge [bend left=30] (n9);
\path (n14) edge [bend left=30] [color = red] (n15);
\path (n8) edge [bend left=40] (n12);

\path (n11) edge [bend right=40] [color = red] (n15);

\path (n8) edge [bend left=50][color = red] (n15);
\path (n1) edge [bend left=20] (n12);
\path (n1) edge [bend left=30] [color = red] (n15);
\path (n3) edge [bend right=20] (n12);
\path (n3) edge [bend right=30][color = red] (n15);
\path (n11) edge [bend right=40] (n12);
\path (n14) edge [color=red] (n12);
\path (n12) edge [color=red] (n16);
\path (n17) edge [color=red] (n15);
\path (n17) edge [color=red] (n16);

 \foreach \from/\to in {n4/n5,n5/n1,n1/n2,n2/n3,n3/n4,n2/n7,n7/n8,n8/n6,n3/n6,n1/n6,n6/n10,n10/n11,n11/n9,n1/n9,n3/n9,n9/n13,n13/n14,n14/n12}
 \draw (\from) -- (\to);
\end{tikzpicture}
\vspace{-6mm}

\caption{Graph $\Gamma_{17}$.}\label{fig:gamma17}
\end{figure}

This graph has in total $5k+\frac{(k-2)(k-1)}{2}$ edges.
For example, below we display graphs~$\Gamma_{14}$ and~$\Gamma_{17}$ to illustrate how to inductively construct $\Gamma_{3k+2}$ from $\Gamma_{3k-1}$ (additional edges and vertices are denoted in red), see Figures~\ref{fig:gamma14} and~\ref{fig:gamma17}.

We first compute the Hilbert series of $R_{\Gamma_{3k+2}}$.
\begin{Theorem}
For $k \geq 1$, we have
\begin{equation} \label{catherine}
H_t(R_{\Gamma_{3k+2}})=\frac{(1+t)^{k-1}(1+(2+k)t+t^2)}{(1-t)^{k+1}}.
\end{equation}
In particular, the Krull dimension of $R_{\Gamma_{3k+2}}$ is $k+1$.
\end{Theorem}
\begin{proof} We require some auxiliary graphs and their edge ideals.
We denote by $\bar{\Gamma}_{3k+2}$ the graph obtained from $\Gamma_{3 k+2}$ by adding a path of length two to the $(3k+2)$-th vertex, i.e.,
$\bar{\Gamma}_{3k+2}$ is obtained from $\Gamma_{3k+2}$ by adjoining a path of length two at the terminal vertex.
For example, in Figure~\ref{fig:gamma17}, $\bar{\Gamma}_{14}$ is the subgraph induced by the vertices $1,\dots,16$. Here a ``path of length two'' means two consecutive edges joining the old terminal vertex to two newly added vertices.
As such, this new graph has $3k+4$ vertices.

Next, we choose in Lemma \ref{recursion} the vertex $v$ to be the $(3k+2)$-th vertex, and we obtain a~relation
\[
H_t(R_{\Gamma_{3k+2}})=H_t(R_{\bar{\Gamma}_{3k-1}})+\frac{t(1+t)^k}{(1-t)^{k+1}},
\]
where for the second term we used (several times) that $H_t(R_{\bullet-\bullet})=\frac{1+t}{1-t}$.

We also apply Lemma \ref{recursion} starting from $\bar{\Gamma}_{3k+2}$ where we take again $v$ to be the $(3k+2)$-th node.
This gives a new relation
\[
H_t(R_{\bar{\Gamma}_{3k+2}})=\frac{1+t}{1-t} H_t(R_{\bar{\Gamma}_{3k-1}})+\frac{t(1+t)^k}{(1-t)^{k+2}}.
\]
We claim that \smash{$H_t(R_{\bar{\Gamma}_{3k+2}})=\frac{(1+t)^{k-1}}{(1-t)^{k+2}}(1+(k+3)t+(k+2)t^2)$}. This follows easily by induction using
the second relation and the fact that \smash{$H_t(R_{\bar{\Gamma}_5})=\frac{1+4t+3t^2}{(1-t)^3}$}, which is easily checked.
The~rest follows from the first relation and the formula for $H_t(R_{\bar{\Gamma}_{3k-1}})$.
\end{proof}

Recall that a standard graded ring is Cohen--Macaulay if its depth equals its Krull dimension. It is Gorenstein if, in addition, its canonical module is isomorphic to a graded shift of the ring itself (equivalently, the canonical module is self-dual in the standard graded setting used here).
The graph families studied here were originally selected because
their edge algebras possess these distinguished algebraic properties. In
particular, among the cycle graphs only $C_5$ is Gorenstein~\cite{Edge}, and it is
precisely this graph that gives rise to the remarkably simple graph series
appearing in~\cite{BJM}.

This example motivates the analogous question for the family $\Gamma_{3k+2}$ appearing in Conjecture~\ref{Conjecture2.4}. Based on the palindromic property of the $h$-coefficients in (\ref{catherine}) and analysis for small~$k$, we expect the following.

\begin{Conjecture}\label{Conjecture2.4} The algebra $R_{\Gamma_{3k+2}}$ is Gorenstein for $k \geq 1$.
\end{Conjecture}

We now turn to the study of the graph series associated with $\Gamma_{3k+2}$.
One drawback of working with leafless graphs is that their graph series are more difficult to analyze, and in particular, we lack an analogue of Lemma~\ref{recursion}.
As already suggested in \cite{BJM}, there does not appear to be a~satisfactory closed combinatorial formula for the graph series of the $n$-cycle, except in the cases $n=3$ and $n=5$.
The following result was used in \cite{BJM}.
\begin{Lemma} \label{double}
For any integers $a,b \geq 0$, we have
\[
\sum_{n,m \geq 0} \frac{q^{mn+m(a+1)+n(b+1)}}{(q)_{m}(q)_n}=\frac{1}{\bigl(q^{b+1}\bigr)_{a+1}\bigl(q^{a+1}\bigr)_{\infty}}
=\frac{1}{\bigl(q^{a+1}\bigr)_{b+1}\bigl(q^{b+1}\bigr)_{\infty}}.
\]
\end{Lemma}
We also record an elementary identity for $q$-Pochhammer symbols:
\begin{equation} \label{poch}
(q)_{m_1+\cdots + m_k}=(q)_{m_1} \bigl(q^{m_1+1}\bigr)_{m_2} \cdots \bigl(q^{m_1+m_2+\cdots + m_{k-1}}\bigr)_{m_k}.
\end{equation}

To compute $H_{\Gamma_{3k+2}}(q)$ as in Theorem~\ref{main}, we require the following result.
\begin{Theorem} \label{3k2} For $k \geq 1$, we have
\begin{align*}
& H_{\Gamma_{3k+2}}(q) =\\
& \sum_{n_{2},n_{5},\dots,n_{3k+2} \geq 0} \frac{q^{\sum_{i=0}^k n_{3i+2}}}{\bigl(1-q^{n_{3k-1}+n_{3k+2}+1}\bigr)\bigl(1-q^{n_{3k-4}+n_{3k-1}+n_{3k+2}+1}\bigr) \cdots \bigl(1-q^{{1+\sum_{i=0}^k n_{3i+2}}}\bigr)},
\end{align*}
where the summation is over $k+1$ variables.
\end{Theorem}
\begin{proof} The theorem clearly holds for $k=1$, where we recover the result from \cite{BJM} so we assume $k \geq 2$.
Next we group certain summation variables in pairs. We first sum over $n_1$ using Euler's identity, and then we sum over $k$ pairs of variables $(n_{3},n_{4})$, $(n_{6},n_{7})$,\dots, and $(n_{3k},n_{3k+1})$ using
Lemma \ref{double}. Notice that this will result in a summation over $k+1$ variables appearing in the formula
\begin{align*}
H_{\Gamma_{3k+2}}(q)={}&\sum_{n_1,\dots,n_{3k+2} \geq 0} \frac{1}{\prod_{i=1}^k (q)_{n_{3i-1}}} \\
&
\times \frac{q^{n_1(n_2+n_5+\cdots +n_{3k+2}+1)+n_2+n_5+\cdots +n_{3k+2}}}{(q)_{n_1}} \frac{q^{n_3 n_4+n_4(n_5+n_8+\cdots + n_{3k+2}+1)+n_3(n_2+1)}}{(q)_{n_3}(q)_{n_4}}
\\
&
\times \cdots \times \frac{q^{n_{3k} n_{3k+1}+n_{3k+1}(n_{3k+2}+1)+n_{3k}(n_{3k-1}+1)}}{(q)_{n_{3k}}(q)_{n_{3k+1}}}.
\end{align*}
From this, we easily get using Lemma \ref{double}
\begin{align*}
& \sum_{n_2,n_5,\dots,n_{3k+2} \geq 0} \frac{(q)_{\sum_{k=0}^{k} n_{3k+2}} q^{\sum_{k=0}^{k} n_{3k+2}}}{(q)_\infty^{k+1} (q)_{n_{3k+2}}\bigl(q^{n_{3k+2}+1}\bigr)_{n_{3k-1}+1} \bigl(q^{n_{3k+2}+n_{3k-1}+1}\bigr)_{n_{3k-4}+1}} \\
&\qquad\quad {}\times \cdots \times \frac{1}{\bigl(q^{\sum_{i=2}^k n_{3i+2}+1}\bigr)_{n_5+1} \bigl(q^{\sum_{i=1}^k n_{3i+2}+1}\bigr)_{n_2+1}}
\\
&\qquad=\sum_{n_2,n_5,\dots,n_{3k+2} \geq 0} \frac{q^{\sum_{i=0}^k n_{3i+2}}}{(q)_\infty^{k+1} \bigl(1-q^{n_{3k+2}+n_{3k-1}+1}\bigr)\bigl(1-q^{n_{3k+2}+n_{3k-1}+n_{3k-4}+1}\bigr)}\\
&\qquad\quad {}
\times \cdots \times \frac{1}{\bigl(1-q^{\sum_{i=1}^k n_{3i+2}+1}\bigr) \bigl(1-q^{\sum_{i=0}^k n_{3i+2}+1}\bigr)},
\end{align*}
where in the last line we used the elementary identity which follows directly from (\ref{poch}):
\begin{align*}
& \frac{(q)_{\sum_{i=0}^k n_{3i+2}}}{(q)_{n_{3k+2}}\bigl(q^{n_{3k+2}+1}\bigr)_{n_{3k-1}+1} \bigl(q^{n_{3k+2}+n_{3k-1}+1}\bigr)_{n_{3k-4}+1} \cdots \bigl(q^{\sum_{i=1}^k n_{3i+2}+1}\bigr)_{n_2+1}} \\
& \qquad =\frac{(q)_{n_{3k+2}}\bigl(q^{n_{3k+2}+1}\bigr)_{n_{3k-1}} \bigl(q^{n_{3k+2}+n_{3k-1}+1}\bigr)_{n_{3k-4}} \cdots \bigl(q^{\sum_{i=1} n_{3i+2}+1}\bigr)_{n_2}}{(q)_{n_{3k+2}}\bigl(q^{n_{3k+2}+1}\bigr)_{n_{3k-1}+1} \bigl(q^{n_{3k+2}+n_{3k-1}+1}\bigr)_{n_{3k-4}+1} \cdots \bigl(q^{\sum_{i=1} n_{3i+2}+1}\bigr)_{n_2+1}} \\
& \qquad =\frac{1}{\bigl(1-q^{n_{3k+2}+n_{3k-1}+1}\bigr) \bigl(1-q^{n_{3k+2}+n_{3k-1}+n_{3k-4}+1}\bigr) \cdots \bigl(1-q^{\sum_{i=0}^k n_{3i+2}+1}\bigr)}.\tag*{\qed}
\end{align*}\renewcommand{\qed}{}
\end{proof}

To prove Theorem~\ref{3k2}, we first require an auxiliary result.
\begin{Lemma} \label{qz} We have
\[
\sum_{n_1 \geq n_2 \geq \cdots \geq n_k \geq 1} \frac{q^{n_1} x^{n_k}}{(1-q^{n_1}) \cdots (1-q^{n_k})}=\sum_{n_1 \geq n_2 \geq \cdots \geq n_k \geq 1} \frac{x q^{n_1}}{(1-x q^{n_1})(1-q^{n_2}) \cdots (1-q^{n_k})}.
\]
\end{Lemma}
\begin{proof} We compare the coefficients of $x^{m}$, $m \geq 1$, on both sides. On the left-hand side, we clearly have
$\sum_{n_1 \geq n_2 \geq \cdots \geq n_{k-1} \geq m} \frac{q^{n_1}}{(1-q^{n_1}) \cdots (1-q^{n_{k-1}})(1-q^{m})}$. On the right-hand side, we
get (using $\frac{xq^{n_1}}{1-x q^{n_1}}=\sum_{\ell \geq 1} x^\ell q^{n_1 \ell}$)
\begin{align*}
 \sum_{n_1 \geq n_2 \geq \cdots \geq n_k \geq 1} \frac{q^{n_1 m}}{(1-q^{n_2}) \cdots (1-q^{n_k})}
 =\sum_{ n_2 \geq \cdots \geq n_k \geq 1} \frac{q^{n_2 m}}{(1-q^{n_2}) \cdots (1-q^{n_k})(1-q^m)}.
\end{align*}
Thus it remains to show
\[
\sum_{n_1 \geq n_2 \geq \cdots \geq n_{k} \geq m} \frac{q^{n_1}}{(1-q^{n_1}) \cdots (1-q^{n_{k}})}=\sum_{ n_1 \geq \cdots \geq n_k \geq 1} \frac{q^{n_1 m}}{(1-q^{n_1}) \cdots (1-q^{n_k})}.
\]
This certainly holds for $m=1$ (for all $k$). Next we write
\begin{align*}
& \sum_{ n_1 \geq \cdots \geq n_k \geq 1} \frac{q^{n_1 (m+1)}}{(1-q^{n_1}) \cdots (1-q^{n_k})}=\sum_{ n_1 \geq \cdots \geq n_k \geq 1} \left(\frac{1}{1-q^{n_1}}-1 \right) \frac{q^{n_1 m}}{(1-q^{n_2}) \cdots (1-q^{n_k})} \\
& \qquad = \sum_{ n_1 \geq \cdots \geq n_k \geq 1} \frac{q^{n_1 m}}{(1-q^{n_1}) \cdots (1-q^{n_k})} - \sum_{ n_2 \geq \cdots \geq n_{k} \geq 1} \sum_{n_1 \geq n_2}
\frac{q^{n_1 m}}{(1-q^{n_1}) \cdots (1-q^{n_k})} \\
& \qquad = \sum_{ n_1 \geq \cdots \geq n_k \geq 1} \frac{q^{n_1 m}}{(1-q^{n_1}) \cdots (1-q^{n_k})} - \sum_{ n_2 \geq \cdots \geq n_{k} \geq 1}
\frac{q^{n_2 m}}{(1-q^{n_2}) \cdots (1-q^{n_k})(1-q^m)}.
\end{align*}
Also observe
\begin{align*}
& \sum_{n_1 \geq n_2 \geq \cdots \geq n_{k} \geq m+1} \frac{q^{n_1}}{(1-q^{n_1}) \cdots (1-q^{n_{k}})} \\
& \qquad =\Biggl( \sum_{n_1 \geq n_2 \geq \cdots \geq n_{k} \geq m} - \sum_{n_1 \geq n_2 \geq \cdots \geq n_{k-1} \geq m \atop n_{k}= m} \Biggr) \frac{q^{n_1}}{(1-q^{n_1}) \cdots (1-q^{n_{k}})} \\
&\qquad =\sum_{n_1 \geq n_2 \geq \cdots \geq n_{k} \geq m} \frac{q^{n_1}}{(1-q^{n_1}) \cdots (1-q^{n_{k}})} \\
& \qquad\quad-\sum_{n_1 \geq n_2 \geq \cdots \geq n_{k-1} \geq m} \frac{q^{n_1}}{(1-q^{n_1}) \cdots (1-q^{n_{k-1}})(1-q^m)}.
\end{align*}
The proof now follows by induction on $m$.
\end{proof}

Taking the first derivative $\frac{{\rm d}}{{\rm d}x}$ and then letting $x=1$ in Lemma~\ref{qz} gives
\begin{equation*}
\sum_{n_1 \geq n_2 \geq \cdots \geq n_k \geq 1} \frac{n_k q^{n_1}}{(1-q^{n_1}) \cdots (1-q^{n_k})}=\sum_{n_1 \geq \cdots \geq n_k \geq 1} \frac{q^{n_1}}{(1-q^{n_1})^2(1-q^{n_2})\cdots (1-q^{n_k})}.
\end{equation*}
To compare with the graph series of $\Gamma_{3k+2}$, the left-hand side above can be written as
\begin{align}
& \sum_{n_1 \geq n_2 \geq \cdots \geq n_k \geq 1} \frac{n_k q^{n_1}}{(1-q^{n_1}) \cdots (1-q^{n_k})} \nonumber\\
& = \sum_{{\bf m} \in \mathbb{N}^{k+1}} \frac{q^{m_1+\cdots + m_{k+1}+1}}{(1-q^{m_2+m_3+1}) \cdots (1-q^{m_2+m_3+ \cdots m_{k+1}+1})(1-q^{m_1+m_2+m_3+ \cdots m_{k+1}+1})}.\label{nk}
\end{align}
The last relation easily follows by letting
\begin{align*}
&n_{k} =m_2+m_3+1, \\
& n_{k-1} =m_2+m_3+m_4+1, \\
 & \cdots \\
& n_1=m_1+m_2+m_3+ \cdots +m_{k+1}+1.
\end{align*}
 The linear factor in the numerator in (\ref{nk}) is the factor $n_k$, which after the change of variables equals the number of possible splittings of the relevant summation variable into $m_2$ and $m_3$.

Following the notation from the introduction, we have $\mathfrak{z}_q(m)=\mathfrak{z}^*_q(m)=\sum_{n \geq 1} \frac{q^n}{(1-q^n)^m}$ and
\[
\mathfrak{z}_q^*(2,1,\dots,1)=\sum_{n_1 \geq n_2 \geq \cdots \geq n_k \geq 1} \frac{q^{n_1}}{(1-q^{n_1})^2 (1-q^{n_2}) \cdots (1-q^{n_k})}.
\]
\begin{Lemma}We have
\[
\mathfrak{z}_q^*(\underbrace{2,1,\dots,1}_{\text{$k$-tuple}})=\sum_{n \geq 1} \frac{n {n+k-2 \choose k-1} q^n}{1-q^n}.
\]
\end{Lemma}
\begin{proof}
This, for instance, follows directly as a special case of the main theorem in \cite{Ohno} (cyclic formula)
\[
\mathfrak{z}_q^*(\underbrace{2,1,\dots,1}_{\text{$k$-tuple}})=k \cdot \mathfrak{z}_q^*(k+1)-(k-1) \cdot \mathfrak{z}_q^*(k)
\]
and properties of binomial coefficients.
\end{proof}

This concludes the proof of Theorem~\ref{main}\,(i).

\begin{Remark}One can show that the adjacency matrix of $\Gamma_{3k+2}$ has determinant $(-1)^{k+1}{(k+1)}$.
Consequently, the graph series associated with $\Gamma_{3k+2}$ may be interpreted as the character of the principal subspace $W_{\Gamma_{3k+2}}$ corresponding to a non-degenerate lattice $L$ of rank $3k+2$; see~\cite{MP}.
It would be interesting to investigate $q$-series identities arising from distinguished $W_{\Gamma_{3k+2}}$-modules associated with irreducible $V_L$-modules.
For the case $k=1$, such identities were studied in~\cite{BJM}.
\end{Remark}

\subsection[Graph series of star graphs T\_\{2k+2\}]{Graph series of star graphs $\boldsymbol{T_{2k+2}}$}\label{section2.3}

Recall that an $\ell$-star graph is a graph (more precisely, a tree) consisting of a central vertex of valency~$\ell$ together with $\ell$ legs of arbitrary length.
In particular, every $\ell$-star graph has exactly~$\ell$ leaves.
An $\ell$-star graph is called \emph{simple} if all of its legs have length one.
The graph series associated with simple $\ell$-star graphs were studied in~\cite{BJM}.

\begin{figure}[!ht]\centering
\begin{tikzpicture}[scale=1.3]
 \node (n1) at (9,4) {1};
 \node (n3) at (7,4) {3};
 \node (n2) at (8,4) {2};
 \node (n4) at (10,4) {4};
 \node (n5) at (11,4) {5};
 \node (n10) at (9,5.2) {10};
 \node (n7) at (7,5.2) {7};
 \node (n6) at (8,4.6) {6};
 \node (n8) at (10,4.6) {8};
 \node (n9) at (11,5.2) {9};
 \foreach \from/\to in {n3/n2,n2/n1,n1/n4,n4/n5,n1/n6,n6/n7,n1/n8,n8/n9,n1/n10}
 \draw (\from) -- (\to);
\end{tikzpicture}\vspace{-1mm}

\caption{$T_{10}$ graph.} \label{T10}
\end{figure}
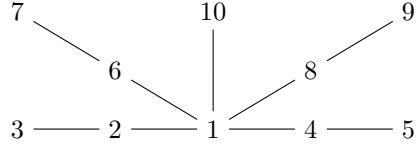

Next, we define a family of graphs $T_{2k+2}$, for $k \geq 2$, each having $2k+2$ vertices.
The graph $T_{2k+2}$ is a $(k+1)$-star graph with $k$ legs of length two and one leg of length one.
We label the central vertex by~$1$.
The leaf vertices are labeled by $3,5,\dots,2k+1$ and $2k+2$, where the vertex $2k+2$ corresponds to the short leg.
The vertices $2,4,6,\dots,2k$ all have valency two.
For $k=2$, the graph $T_6$ coincides with the Dynkin diagram of type $E_6$, while for $k=4$ the corresponding graph is depicted in Figure~\ref{T10}.

The Hilbert series of $R_{T_{2k+2}}$ is easily computed giving \smash{$H_t(R_{T_{2k+2}})=\frac{t+(1+t)^k}{(1-t)^{k+1}}$}.

To compute its graph series, and finish proving Theorem~\ref{main}, we only need basic manipulations with $q$-series.
\begin{Proposition} We have
\[
H_{T_{2k+2}}(q)=q^{-1}\frac{\sum_{n \geq 1} \dfrac{{n+k-2 \choose k-1}q^n}{1-q^n}}{(q)_\infty^{k+1}}.
\]
\end{Proposition}
\begin{proof} By definition, the graph series of $T_{2k+2}$ is
\[
H_{T_{2k+2}}(q)
=
\sum_{n_1,\dots,n_{2k+2}\geq0}
\frac{
q^{
\sum_{i=1}^{2k+2}n_i
+n_1n_{2k+2}
+\sum_{j=1}^{k} (n_1n_{2j}+n_{2j}n_{2j+1} )
}
}{
(q)_{n_1}(q)_{n_2}\cdots(q)_{n_{2k+2}}
}.
\]
First, we sum over the variables $n_3,n_5,n_7,\dots,n_{2k+1}. $ For each of these variables, the summand has the form $ \frac{q^{r(n_i+1)}}{(q)_r},$ where $r$ denotes the variable being summed and $n_i$ is the adjacent remaining variable. By the $q$-binomial identity \[ \sum_{r\geq0}\frac{z^r}{(q)_r} = \frac{1}{(z;q)_\infty}, \] we obtain \[ \sum_{r\geq0}\frac{q^{r(n_i+1)}}{(q)_r} = \frac{1}{(q^{n_i+1};q)_\infty} = \frac{(q)_{n_i}}{(q)_\infty}. \] The numerator $(q)_{n_i}$ cancels the corresponding factor already present in the denominator of the original summand. Since there are $k$ such summations, their total contribution is $\frac{1}{(q)_\infty^k}. $ We next sum over the remaining intermediate variables $ n_2,n_6,n_8,\dots,n_{2k}.$ After the preceding cancellations, each of these variables occurs only through a geometric factor $q^{r(n_1+1)}$. Hence $ \sum_{r\geq0}q^{r(n_1+1)} = \frac{1}{1-q^{n_1+1}}. $ There are $k$ such geometric sums, and therefore together they contribute $ \frac{1}{(1-q^{n_1+1})^k}. $

To prove the assertion, we observe
\begin{align*}
\sum_{n \geq 1} \dfrac{q^n}{(1-q^n)^k}& =\sum_{n \geq 1} \sum_{\ell \geq 0}(-1)^\ell {-k \choose \ell } q^{n (\ell+1)}\\
& = \sum_{\ell \geq 1, n \geq 1} {k+\ell-2 \choose \ell-1} q^{n \ell}=\sum_{\ell \geq 1} \frac{{k+\ell-2 \choose k-1} q^\ell}{1-q^\ell}.\tag*{\qed}
\end{align*}\renewcommand{\qed}{}
\end{proof}

The edge algebra of $T_{2k+2}$ is not always Gorenstein (e.g., $T_{10}$). However, we still have the following.

\begin{Proposition} \label{CM} The edge algebra $R_{T_{2k+2}}$ is Cohen--Macaulay for every $k \geq 1$.
\end{Proposition}
\begin{proof}
As already pointed out, $T_{2k+2}$ is a bipartite graph with the bipartition
\[
X = \{1, 3, \dots, 2k+1\} \qquad \text{and} \qquad Y = \{2, 4, \dots, 2k+2\}.
\]
We relabel the vertices as
\[
x_i := 2i - 1, \quad 1 \leq i \leq k+1,
\qquad \text{and} \qquad
y_i := 2k + 4 - 2i, \quad 1 \leq i \leq k+1.
\]
With this notation, the following clearly hold:
\begin{enumerate}\itemsep=0pt
 \item[(1)] $x_i \sim y_i$ for every $i$,
 \item[(2)] if $x_i \sim y_j$, then $i \leq j$,
 \item[(3)] for $i < j < k$, the condition $x_i \sim y_j$ and $x_j \sim y_k$ implies $x_i \sim y_k$. However, this never occurs since $x_i \sim y_j$ and $x_j \sim y_k$ does not happen simultaneously.
\end{enumerate}
These conditions (1)--(3) are necessary and sufficient for $R_\Gamma$ (when $\Gamma$ is bipartite) to be Cohen--Macaulay \cite[Theorem~3.4]{CM}.
\end{proof}

\subsection[Graph series and multiple q-zeta series]{Graph series and multiple $\boldsymbol{q}$-zeta series}\label{section2.4}

Now we are ready to explain how to construct all multiple $q$-zeta values in the framework of graph series and thus of vertex algebras. As already discussed in the introduction, we
are interested
in the $q$-MZV model given by
\[
{\mathfrak{z}}_q^*(a_1,\dots,a_k)=\sum_{n_1 \geq n_2 \geq \cdots \geq n_k \geq 1} \frac{q^{n_1}}{(1-q^{n_1})^{a_1} \cdots (1-q^{n_k})^{a_k}}.
\]
Consider a disconnected graph with $2k$ vertices $\{1,\dots,k,1',\dots,k' \}$ and $k$ edges:
\begin{center}
\begin{tikzpicture}[scale=1.2]
 \node (n1) at (1,1) {1};
 \node (n11) at (1.5,1) {1$'$};
 \node (n2) at (2,1) {2};
 \node (n12) at (2.5,1) {2$'$};
 \node (n3) at (3,1) {3};
 \node (n13) at (3.5,1) {3$'$};
 \node (n4) at (4,1) {4};
 \node (n14) at (4.5,1) {4$'$};
 \node (n5) at (5,1) {5};
 \node (n15) at (5.5,1) {5$'$};
 \node (n6) at (6,1) {...};
 \node (n7) at (7,1) {...};
 \node (n8) at (8,1) {...};
 \node (n9) at (9,1) {...};
 \node (n20) at (10.5,1) {$k'$};
 \node (n10) at (10,1) {$k$};

 \foreach \from/\to in {n1/n11,n2/n12,n3/n13,n4/n14,n5/n15,n10/n20}
 \draw (\from) -- (\to);
\end{tikzpicture}
\end{center}
We next construct a graph, denoted by $Z_{a_1,\dots,a_k}$, which will contain the above graph as a subgraph. In this new graph, the nodes denoted by $1',\dots,k'$ are leaves (so they are only adjacent to $1,2,\dots,k$, respectively).

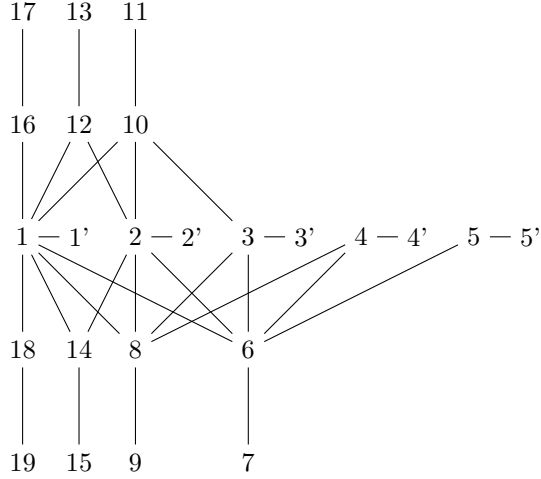
\begin{figure}[!ht]\centering
\begin{tikzpicture}[scale=1.5]
 \node (n1) at (1,1) {1};
 \node (n21) at (1.5,1) {1$'$};
 \node (n2) at (2,1) {2};
 \node (n22) at (2.5,1) {2$'$};
 \node (n3) at (3,1) {3};
 \node (n23) at (3.5,1) {3$'$};
 \node (n4) at (4,1) {4};
 \node (n24) at (4.5,1) {4$'$};
 \node (n5) at (5,1) {5};
 \node (n25) at (5.5,1) {5$'$};
 \node (n6) at (3,0) {6};
 \node (n7) at (3,-1) {7};
 \node (n8) at (2,0) {8};
 \node (n9) at (2,-1) {9};
 \node (n10) at (2,2) {10};
 \node (n11) at (2,3) {11};
 \node (n12) at (1.5,2) {12};
 \node (n13) at (1.5,3) {13};
 \node (n14) at (1.5,0) {14};
 \node (n15) at (1.5,-1) {15};
 \node (n16) at (1,2) {16};
 \node (n17) at (1,3) {17};
 \node (n18) at (1,0) {18};
 \node (n19) at (1,-1) {19};

 \foreach \from/\to in {n1/n18,n1/n14,n1/n8,n1/n6,n1/n16,n1/n12,n1/n10,n2/n12,n2/n10,n2/n14,n2/n8,n2/n6,n3/n10,n3/n8,n3/n6,n4/n8,n4/n6,n5/n6,n18/n19,n14/n15,n8/n9,n6/n7,n16/n17,n12/n13,n10/n11,n1/n21,n2/n22,n3/n23,n4/n24,n5/n25}
 \draw (\from) -- (\to);
\end{tikzpicture}
\caption{Graph $Z_{1,1,1,2,2}$.}\label{figure5}
\end{figure}

To construct our graph, we attach various simple star graphs to the underlying graph with~$2k$ vertices. For better transparency, in parallel, we consider the graph
$Z_{1,1,1,2,2}$ displayed on Figure~\ref{figure5}.

We first attach $a_1$ simple $(k+1)$-star graphs to the $k$ nodes enumerated by $1,2,\dots,k$ (e.g., on Figure~\ref{figure5}, we attach one $6$-star graph with the central node $6$ to the nodes $1$, $2$, $3$, $4$, $5$).
Then we attach $a_2$ simple star graphs with $k$ legs to the $k-1$ nodes enumerated by $1,2,\dots,k-1$ (e.g., on Figure~\ref{figure5}, we attach one $5$-star graph with central node $8$ to the nodes $1$, $2$, $3$, $4$). We continue until we attach $a_k$ $2$-star graphs (i.e., paths of length two) to node~1 (e.g., on Figure~\ref{figure5}, we attach two $2$-star graphs with central vertices $16$ and $18$ to node~$1$).

The Hilbert series of $R_{Z_{a_1, \dots,a_k}}$ can be recursively computed using the relations
\begin{align*}
&H_t(R_{Z_{a_1, \dots,a_k}}) =H_t(R_{Z_{a_1-1, \dots,a_k}})+\frac{t}{(1-t)^{k+1}} \left(\frac{1+t}{1-t}\right)^{a_1+\cdots +a_k-1}, \\
&H_t(R_{Z_{0, \dots,0,a_i,a_{i+1}, \dots, a_k}}) =H_t(R_{Z_{a_i, \dots,a_k}}).
\end{align*}
A closed formula for $H_t(R_{Z_{a_1, \dots,a_k}})$ is somewhat messy so we do not display it here.

For the graph series, we have the following result:
\[
H_{Z_{a_1,\dots,a_k}}(q)=\frac{q^{-1}{\mathfrak{z}}^*_q(a_1,\dots,a_k)}{(q)_\infty^{a_1+\cdots +a_k+k}}.
\]
\begin{proof}[Proof of Theorem~\ref{qmzv}]
Again, the proof is easier to follow by considering an example in~parallel; see Figure~\ref{figure5}.
We first sum over all $(k+1)$-star subgraphs. By definition, there are $a_1$ such subgraphs.
For each subgraph, we start from the leaf node $s$ (for instance, node 7 in Figure~\ref{figure5}) and the central node $t$ (for instance, node 6 in Figure~\ref{figure5}).
This leads to the following summation:
\begin{align*}
\sum_{n_s,n_t \geq 0} \frac{q^{n_s(n_t+1)+n_t(n_1+\cdots + n_k+1)}}{(q)_{n_s}(q)_{n_t}}
=\sum_{n_t \geq 0} \frac{q^{n_t(n_1+\cdots + n_k+1)}}{(q^{n_t+1})_\infty (q)_{n_t} }= \frac{1}{(q)_\infty \bigl(1-q^{n_1+\cdots +n_k+1}\bigr)}.
\end{align*}
Repeating this $a_1$ times (observe that the summation variables corresponding to the relevant subgraphs are independent) results with
\[
\frac{1}{(q)^{a_1}_\infty \bigl(1-q^{n_1+\cdots +n_k+1}\bigr)^{a_1}}.
\]
We repeat the same procedure for $k$-star graphs attached to $1,2,\dots,k-1$. These subgraphs contribute with
\[
\frac{1}{(q)^{a_2}_\infty \bigl(1-q^{n_1+\cdots +n_{k-1}+1}\bigr)^{a_2}}.
\]
We continue until we sum over paths of length two (e.g., nodes 17 and 16) attached to the node $1$. Here the contribution is
\[
\frac{1}{(q)^{a_k}_\infty \bigl(1-q^{n_1+1}\bigl)^{a_k}}.
\]
We also have to sum over the variables that correspond to the ``prime'' nodes:
\[
\sum_{ n_1',\dots,n_k' \geq 0} \frac{q^{n_1'(n_1+1)+\cdots + n_k'(n_k+1)}}{(q)_{n_1'} \cdots (q)_{n_k'}}=\frac{1}{\bigl(q^{n_1+1}\bigr)_\infty \cdots \bigl(q^{n_k+1}\bigr)_\infty}.
\]
Combining everything with the final $k$-fold summation over $n_1,\dots,n_k$, we quickly get
\[
\frac{1}{(q)_\infty^{\sum_{i=1}^k a_i+k}} \sum_{n_1,\dots,n_k \geq 0} \frac{q^{n_1+ \cdots + n_k}}{\bigl(1-q^{\sum_{i=1}^{k} n_i+1}\bigr)^{a_1} \bigl(1-q^{\sum_{i=1}^{k-1} n_i+1}\bigr)^{a_2} \cdots \bigl(1-q^{n_1+1}\bigr)^{a_k}}.
\]
Setting $n_1+n_2+\cdots + n_k+1 \to m_1$, $n_1+n_2+\cdots + n_k+1 \to m_2,\dots,n_1+1 \to m_k$ proves the claim.
\end{proof}

\begin{Remark}
Observe that $Z_{a_1,\dots,a_k}$ is a connected bipartite graph. Using exactly the same type of argument as in the proof Proposition~\ref{CM},
it follows that the rings $R_{Z_{a_1,\dots,a_k}}$ are Cohen--Macaulay.
\end{Remark}

\subsection{Algebra of graphs series}\label{section2.5}
The set of all graph series of simple graphs is closed under the multiplication
\begin{equation} \label{mult}
H_{\Gamma_1}(q) \cdot H_{\Gamma_2}(q)=H_{\Gamma_1 \oplus \Gamma_2}(q).
\end{equation}
Under this operation we form a $\mathbb{Q}$-algebra
\[
\mathcal{G}=\mathbb{Q} \langle H_{\Gamma}(q) \mid \Gamma \in {\rm SGraphs}/{\sim} \rangle,
\]
taking the $\mathbb{Q}$-span over all isomorphism classes of simple graphs.
This algebra is filtered,
\[
\mathcal{G}=\bigcup_{n \geq 0} \mathcal{G}_n, \qquad \mathcal{G}_n \cdot \mathcal{G}_m \subset \mathcal{G}_{m+n},
\]
where $\mathcal{G}_n$ is spanned by all $H_{\Gamma}(q)$ where $\Gamma$ runs through the set of (non-isomorphic)
graphs with $\leq n$ vertices, and $\mathcal{G}_0=\mathbb{Q}$.
As we already indicated, apart from the obvious relations (\ref{mult}), new relations start to emerge for graphs with six and more vertices.
Moreover, this $\mathbb{Q}$-algebra is not graded, as is visibly seen from the (non-homogeneous) relation
\[
H_{\Gamma_8 \oplus pt}=2 H_{T_8}(q)-H_{T_6 \oplus pt}(q),
\]
a consequence of Theorem~\ref{main}. Let ${\rm gr}(\mathcal{G})$ denote the associated graded algebra of $\mathcal{G}$.
Summarizing the previous discussion:
\begin{alignat*}{4}
& \dim ({\rm gr}(\mathcal{G})_1)=1, \qquad&& \dim ({\rm gr}(\mathcal{G})_2)=2, \qquad && \dim ({\rm gr}(\mathcal{G})_3)=4,&\\
& \dim ({\rm gr}(\mathcal{G})_4)=11,\qquad&& \dim ({\rm gr}(\mathcal{G})_5)=34, \qquad && &
\end{alignat*}
all coinciding with the number of non-isomorphic graphs on $1 \!\leq\! n\! \leq\! 5$ vertices, but $\dim ({\rm gr}(\mathcal{G})_6)<156$, the number of non-isomorphic graphs on six vertices.
Most graph series are unfamiliar $q$-series even for graphs with a small number of vertices.
For instance, among the 21 connected graphs with five vertices, their graph series $H_\Gamma(q)$, modulo multiplicative powers of the Euler product, do not appear in the OEIS database~\cite{OEIS}, except for $H_{C_5}(q)$ and $H_{D_5}(q)$, which were studied in~\cite{BJM}.

\subsection{Quasi-modularity and graph series}\label{section2.6}
We end this part with some comments on quasi-modularity of $q$-MZVs.
This feature is important and directly related to conjectures on graded dimensions
of the algebra of $q$-MZVs. For classical multiple zeta values, we
have a well-known formula \cite{Hoffman}:
\[
\sum_{ \sigma \in S_n} \zeta(2k_{\sigma(1)},\dots,2k_{\sigma(n)}) \in \mathbb{Q}[\pi^2].
\]
More precisely, the left-hand side can be expressed as a sum of products of the ordinary zeta values $\zeta(2k)$. For $q$-MZVs, this type of formula is formulated in terms of $\zeta_q(2k)$
or, equivalently, in terms of quasi-modular forms. Either way, \cite{Bradley} shows that
\[
\mathcal{\zeta}_q(2k,2k,\dots,2k)
\]
can be expressed as a polynomial in $q$-zeta values $\zeta_q(m)$. Even though \cite{AR}
\[
\mathcal{\zeta}_q(2,2,\dots,2) \in \mathcal{QM},
\]
 in general
\[
\mathcal{\zeta}_q(2k,2k,\dots,2k) \notin \mathcal{QM}.
\]

\section[Multiple q-zeta values, quasi-modular forms, and characters]{Multiple $\boldsymbol{q}$-zeta values, quasi-modular forms, and characters}\label{section3}

This section contains the second part of the paper. Throughout, we let
$\tau\in\mathbb H$, $z\in\mathbb C$, $q={\rm e}^{2\pi {\rm i}\tau}$, and
$\zeta={\rm e}^{2\pi {\rm i} z}$.
Unless stated otherwise, all identities involving Fourier series in $\zeta$
are interpreted in the annulus specified in the statement, which determines
the corresponding geometric-series expansion of each denominator. We begin by
recalling the elliptic and quasi-Jacobi functions needed for our constant-term
computations and then apply these techniques to characters of genus-zero
class-$\mathcal{S}$ vertex operator algebras.

\subsection{Elliptic functions, theta functions, and quasi-modular forms}\label{section3.1}
We present some basic (and mostly well-known) results on elliptic functions, (quasi)-Jacobi forms, and their Fourier coefficients. We closely follow \cite{Lang}.

The Jacobi theta function is given by
\[
\vartheta(z;\tau)=\sum_{n \in \mathbb{Z}+\frac12} {\rm e}^{\pi {\rm i} n^2 \tau+2 \pi {\rm i} n(z+\frac12)}.
\]
We also define Eisenstein series with their $q$-expansion, for integers $k\geq1$ (and with the usual convention that $E_2$ is quasimodular)
\[
G_{2k}(\tau)=2 \zeta(2k)+\frac{(2\pi {\rm i})^{2k}}{(2k-1)!} \sum_{n \geq 1} \frac{n^{2k-1}q^n}{1-q^n}, \qquad q={\rm e}^{2 \pi {\rm i} \tau},
\]
and slightly normalized versions
\begin{equation} \label{eis}
E_{2k}(\tau)=1-\frac{4k}{B_{2k}}\sum_{n \geq 1} \frac{n^{2k-1} q^n}{1-q^n},
\end{equation}
where $B_{2k}$ is the Bernoulli number. Recall that $\mathcal{QM}$ us the ring of quasi-modular forms for the full modular group discussed earlier. A famous result of S. Ramanujan is that
\[
\left(q \frac{{\rm d}}{{\rm d}q} \right)\colon \ \mathcal{QM} \rightarrow \mathcal{QM},
\]
so $\mathcal{QM}$ is a differential ring. This well-known fact is employed throughout this section and in later parts.

The Weierstrass $\zeta$-function\footnote{Not to be confused with the Fourier coefficient $\zeta={\rm e}^{2 \pi {\rm i} z}$ used below.} and the $\wp$-function are defined as
\begin{align*}
&\zeta(z;\tau) =\frac1z+\sum_{(n,m) \neq (0,0)} \left(\frac{1}{z-m \tau-n}+\frac{1}{m\tau+n}+\frac{z}{(m\tau+n)^2} \right),\\
&\wp(z; \tau) =\frac{1}{z^2}+\sum_{(n,m) \neq (0,0)} \frac{1}{(z-m-n\tau)^2}-\frac{1}{(m+n \tau)^2}, \\
&\wp'(z;\tau)=-\frac{1}{2} \frac{{\rm d}}{{\rm d}z} \wp(z;\tau).
\end{align*}
Their Laurent expansions around $z=0$ are then given by
\begin{align*}
&\zeta(z;\tau) =\frac{1}{z}-\sum_{n \geq 1} G_{2n+2}(\tau)z^{2n+1}, \\
&\wp(z;\tau)=\frac{1}{z^2} + \sum_{n \geq 1} (2n+1)G_{2n+2}(\tau) z^{2n}, \\
&\wp'(z;\tau)=\frac{1}{z^3} - \sum_{n \geq 1} n(2n+1) G_{2n+2}(\tau) z^{2n-1}.
\end{align*}
For concrete computation, it is convenient to use the so-called $P$-functions or ``propagators'' (widely used for computations of genus one correlators in CFT and vertex operator algebra~\cite{Zhu}):
\[
P_k(\zeta;q)=\frac{(2 \pi {\rm i})^k}{(k-1)!} \sum_{n \geq 1} \left( \frac{n^{k-1} \zeta^{n}q^n}{1-q^n}+ (-1)^k \frac{n^{k-1} \zeta^{-n}}{1-q^n} \right).
\]
Then inside $|q|<\bigl|\zeta^{-1}\bigr|<1$ we have
\begin{align}
\label{rho1} &P_1\bigl({\rm e}^{2 \pi {\rm i} z},\tau\bigr)=-\zeta(z;\tau)+G_2(\tau)z -\pi {\rm i}, \\
\label{rho2} &P_2\bigl({\rm e}^{2 \pi {\rm i} z}, \tau\bigr)=\wp(z,\tau)+G_2(\tau), \\
\label{rhok} &P_{3}\bigl({\rm e}^{2 \pi {\rm i} z}, \tau\bigr)=-\wp'(z,\tau).
\end{align}
Similarly, $P_k\bigl({\rm e}^{2 \pi {\rm i} z};\tau\bigr)$, $k \geq 4$, is proportional to the $(k-2)$-th $z$-derivative of $\wp(z;\tau)$.
The second and the third relations (\ref{rho2})--(\ref{rhok}) are equivalent to the following:
\begin{align*} 
& \wp(z,\tau)=-G_2(\tau)+(2 \pi {\rm i})^2 \sum_{n \geq 1} \left( \frac{n q^n \zeta^n }{(1-q^n)}+\frac{n \zeta^{-n} }{(1-q^n)} \right),\\
& \wp'(z,\tau)=-(2 \pi {\rm i})^3 \sum_{n \geq 1} \left( \frac{n^2 q^n \zeta^n }{(1-q^n)}-\frac{n^2 \zeta^{-n} }{(1-q^n)} \right).
\end{align*}
To see that (\ref{rho1}) holds, we require an identity
\begin{equation} \label{J1}
-\frac12 +\frac{\zeta}{\zeta-1} - \sum_{n \geq 1}\left(\frac{ \zeta^{-1} q^n}{1- \zeta^{-1}q^n}-\frac{q^n \zeta}{1-q^n \zeta}\right)=\frac12+\sum_{n \geq 1} \left(\frac{\zeta^n q^n}{1-q^n}-\frac{\zeta^{-n}}{1-q^n} \right),
\end{equation}
again valid in $|q|<|\zeta^{-1}|<1$. We denote either side of (\ref{J1}) by $J_1(z;\tau)$.
This function does not transform as a Jacobi form, instead we have
\begin{align*}
&J_1(z+ m \tau+n;\tau)= J_1(z;\tau)-m, \\
&J_1 \left(\frac{z}{\tau};\frac{-1}{\tau}\right)=\tau J_1(z;\tau)+z.
\end{align*}
We also require $\sigma(z;\tau)$, which can be defined using the theta function \cite{Lang}
\[
\sigma(z;\tau)=-{\rm exp}\left(\frac{\eta_1 z^2}{2}\right)\frac{\vartheta(z;\tau)}{2 \pi \eta(\tau)^3},
\]
where $\eta_1=G_2(\tau)$ in the normalization of \cite{Lang}; for our purposes, it is a function of $\tau$ only and will not be used elsewhere.

Functions $\wp$, $\vartheta$, and $\sigma$ are related through a fundamental relation~\cite{Lang}:
\begin{align}
\wp(z_1;\tau)-\wp(z_2;\tau)=\frac{\sigma(z_1+z_2;\tau)\sigma(z_1-z_2;\tau)}{\sigma(z_1;\tau)^2\sigma(z_2;\tau)^2} =\frac{\vartheta(z_1+z_2;\tau) \vartheta(z_1-z_2;\tau)}{ 4 \pi^2 \vartheta(z_1;\tau)^2 \vartheta(z_2;\tau)^2 \eta(\tau)^6}.\label{rho-sigma}
\end{align}

Our next result concerns the constant term in the Fourier expansion of elliptic functions with respect to $\zeta = {\rm e}^{2\pi {\rm i} z}$.
This result is well known and has appeared previously in the literature; see, for example, \cite{GM, MP}.
In the present work, however, we adopt the approach of~\cite{BFM}, where all Fourier coefficients of certain index-zero Jacobi forms are investigated, rather than focusing solely on the constant term.
\begin{Proposition} \label{elliptic} Let $f(z;\tau) \in \mathbb{C}[\wp(z;\tau),\wp'(z;\tau)]$.
Then, in the range $ |q|< \bigl|\zeta^{-1}\bigr|<1$, the constant term ${\rm CT}_\zeta f(z;\tau)$ belongs to $\mathcal{QM}$, and the non-constant terms ${\rm CT}_\zeta f(z;\tau)\zeta^{-r}$, $r \neq 0$, are not quasi-modular.
\end{Proposition}
\begin{proof} We may assume that $f=\wp^r (\wp')^s$. Using the Weierstrass equation
\[
(\wp')^2=4\wp^3-g_2(\tau)\wp-g_3(\tau),
\]
where $g_2(\tau)$ and $g_3(\tau)$ are multiples of $E_4(\tau)$ and $E_6(\tau)$, respectively, every higher power of $\wp'$ can be reduced to the cases $s=0$ and $s=1$.
We use the Laurent expansion in the variable $z$ around zero to write
\begin{align*} 
f(z;\tau)=\sum_{0 \leq i < \frac{\ell}{2}} \frac{D_{\ell-2i}(\tau)}{ z^{\ell-2i}}+D_0(z)+O(z),
\end{align*}
where the $z$-coefficients $D_{\ell-2i}(\tau)$ are holomorphic {modular} forms of weight $2i$ that can be explicitly computed using the Laurent expansion of $\wp$
and $\wp'$, and $\ell=2r+3s$. Of course, $D_0(z)$ could be zero depending on the parity of $s$. To compute the $\zeta$-constant term, we proceed as in \cite[Theorem~1.1]{BFM} dealing with a particular Jacobi form of index zero (the case when $m$ is even in~\cite{BFM} is precisely the case of elliptic functions here).

For the nonzero Fourier coefficient $f_r(\tau)$, $r \neq 0$, we compute the contour integral
\[
\int_{\partial P_{z_0}} f(w;\tau) {\rm e}^{-2 \pi {\rm i} r w} \, {\rm d}w,
\]
where
\[
P_{z_0} = \{ z_0 + r \tau + s \mid 0 \leq r,s \leq 1 \},
\]
in two different ways: directly as a contour integral, and by the residue theorem by choosing
\[
z_0 = -\tfrac{1}{2} - \tfrac{\tau}{2},
\]
so that the only pole inside the contour is at $z=0$.
This quickly leads to a formula for the $r$-th coefficient, which is quasi-modular if and only if $r=0$ and not quasi-modular otherwise (see \cite[Theorem~1.1]{BFM}).
More explicitly,
\[
{\rm CT}_\zeta f(z;\tau)=D_0(\tau)+\sum_{i} \lambda_i D_{2i}(\tau)E_{2i}(\tau),
\]
where $\lambda_i$ are some coefficients. This is certainly contained in~$\mathcal{QM}$.
\end{proof}

\subsection{Quasi-Jacobi forms}\label{section3.2}
We are also interested in Fourier coefficients, or more precisely the constant term, of more general functions that include $\tau$-derivatives of elliptic
functions. For that we need the concept of {\em quasi-Jacobi forms} (see \cite{Libgober,Ob, Ittersum} for instance).
We are only concerned with the subspace of meromorphic quasi-Jacobi forms of index $0$ with poles at $L_\tau=\mathbb{Z}\tau+\mathbb{Z}$. These are known to form a subring $\mathcal{QJ}^0:=\mathbb{C}[E_6,E_4,J_1,\wp,\wp']$ inside the ring of quasi-Jacobi forms $\mathcal{QJ}:=\mathbb{C}[F,E_6,E_4,J_1,\wp,\wp']$ ($E_6$ can be omitted from both rings though), where $F$ is a particular quasi-Jacobi form of non-zero index not needed here \cite{Ob}. Besides $J_1(z;\tau)$ mentioned earlier, we also require more general quasi-Jacobi forms ($k \geq 2$):
\[
J_k(z;\tau)=B_k - k \sum_{n, r \geq 1} r^{k-1}\bigl( \zeta^n + (-1)^k \zeta^{-n}\bigr) q^{nr}.
\]
Next we define the polynomials $p_m(x)$, $m \geq 2$ via
\[
\sum_{n \geq 1} n^{m-1} x^n=\frac{x p_{m}(x)}{(1-x)^{m}}.
\]
For instance, $p_2(x)=1$, $p_3(x)=1+x$, etc.
These are known as (slightly shifted) Euler's polynomials. The next result is elementary and essentially follows from the definition of $p_m(x)$.
\begin{Lemma} \label{euler} We have
\begin{itemize}\itemsep=0pt
\item[$(1)$] For every even $s \geq 0$, and $k \geq 1$,
\[
\sum_{r \geq 1} \frac{r^s q^r p_{2k}(q^r)}{(1-q^r)^{2k}} \in \mathcal{QM}.
\]
\item[$(2)$] For every odd $s \geq 1$, and $k \geq 1$,
\[
\sum_{r \geq 1} \frac{r^s q^r (1+q^r) p_{2k}(q^r)}{(1-q^r)^{2k+1}} \in \mathcal{QM}.
\]
\end{itemize}
\end{Lemma}
\begin{proof} For (1), it is enough to observe:
\begin{align*}
 \sum_{r \geq 1} r^j \frac{q^r p_{a+j}(q^r)}{(1-q^r)^{a+j}}= \sum_{r \geq 1,\, n \geq 1}\! r^j n^{a+j-1} q^{n r}= \left(q \frac{{\rm d}}{{\rm d}q} \right)^j \! \sum_{r \geq 1,\, n \geq 1}\! n^{a-1} q^{n r} =\left(q \frac{{\rm d}}{{\rm d}q}\right)^j\! \sum_{n \geq 1} \frac{n^{a-1} q^n}{1-q^n}.
\end{align*}
This implies that, for every $s$ and $k$ as above, we may choose $a \geq 2$ and $j \geq 0$ such that the sum in (1) is a (higher) derivative of an Eisenstein series, or an Eisenstein series with its constant term removed.
Similarly, any series appearing in (2) can be expressed as a linear combination of derivatives of two Eisenstein series.
In either case, the resulting sum lies in $\mathcal{QM}$.
\end{proof}

\begin{Theorem} \label{quasi} Let $f(z;\tau) \in \mathcal{QJ}^0$. Then in the range $|q|<\bigl|\zeta^{-1}\bigr|<1$, ${\rm CT}_\zeta f(z;\tau) \in \mathcal{QM}$.
\end{Theorem}
\begin{proof}
Compared with Proposition~\ref{elliptic}, where $f$ is doubly periodic in $z$ and belongs to $\mathbb{C}[\wp,\wp']$, the only new issue is the contribution of $J_1(z;\tau)$ and thus it suffices to prove the statement for
\begin{align*}
f(z;\tau)=J_1(z;\tau)^j g(z;\tau),
\end{align*}
where $g(z;\tau) \in \mathbb{C}[\wp,\wp']$ is elliptic and $j \geq 1$.
We use induction on $j$. For $j=1$, we have to
analyze
\[
J_1(z;\tau) g(z;\tau),
\]
where $g(z;\tau)$ is elliptic. We will compute the constant term by expanding $g(z;\tau)=\zeta^r g_r(\tau)$, where $g_0(\tau)$ is quasi-modular (possibly zero) and
non-zero Fourier coefficients are given by~\cite{BFM}
\[
g_r(\tau)=\frac{q^{r}}{1-q^r} \sum_{0 \leq \ell < m/2} r^{m-2\ell-1} D_{\ell}(\tau), \qquad r>0,
\]
and
\[
g_r(\tau)=-\frac{1}{1-q^{-r}} \sum_{0 \leq \ell < m/2} r^{m-2\ell-1} D_{\ell}(\tau), \qquad r<0,
\]
where the parity of $m$ depends on the parity of the $z$ powers in the Laurent expansion of $g(z;\tau)$, and $D_{\ell}(\tau)$ are modular forms as in Proposition~\ref{elliptic}.
Multiplying out $J_1(z;\tau)=\frac12+\sum_{n \geq 1} \bigl(\frac{\zeta^n q^n}{1-q^n}-\frac{\zeta^{-n}}{1-q^n} \bigr)$ and $g(z;\tau)$, we quickly find the constant term of $J_1(z;\tau) g(z;\tau)$ to be
\[
\frac{1}{2}g_0(\tau)- \sum_{0 \leq \ell < m/2} D_\ell(\tau) \sum_{r \geq 1} \frac{ r^{m-2\ell-1} q^r}{(1-q^r)^2} - \sum_{0 \leq \ell < m/2} D_\ell(\tau) \sum_{r \geq 1} \frac{ (-r)^{m-2\ell-1} q^r}{(1-q^r)^2}.
\]
If $m$ is even, then the two sums cancel out and we end up with $\frac12 g_0(\tau)$ for the constant term which is certainly in $\mathcal{QM}$ by Proposition~\ref{elliptic}. If $m$ is odd, then we get a contribution from the sum
\[
\sum_{r \geq 1} \frac{ r^{s} q^r}{(1-q^r)^2},
\]
where $s:=m- 2\ell -1 \geq 0$ is even. For $s=0$, this sum equals $\frac{1}{24}(1-E_2(\tau)) \in \mathcal{QM}$. For $s \geq 2$, the sum is the first $\bigl(q \frac{{\rm d}}{{\rm d}q}\bigr)$-derivative of
the Eisenstein series $E_{s-1}(\tau)$ and thus in $\mathcal{QM}$.

For the induction step, we need \cite[formula~(11)]{Ob-Serre} that allows us to express $J_1^\ell$ as a linear combination of $J_\ell$, an elliptic Jacobi form $K_\ell$, and quasi-Jacobi forms involving $J_1^k$, $k < \ell$, for which the induction hypothesis applies.
Thus in order to finish the proof we have to argue that
\[
{\rm CT}_\zeta J_\ell(z;\tau) g(z;\tau) \in \mathcal{QM},
\]
where $g(z;\tau)$ is elliptic. As in the $j=1$ case, we first observe that
\[
\frac{1}{\ell}(B_{\ell}-J_{\ell}(z;\tau))=\sum_{n \geq 1}\left( \frac{q^n p_\ell(q^n) \zeta^n}{(1-q^n)^\ell} + (-1)^\ell \frac{q^n p_\ell(q^n) \zeta^{-n}}{(1-q^n)^\ell} \right),
\]
where $p_\ell$ are as above. We distinguish two cases: $\ell$ is even and $\ell$ is odd.
Assume first that $\ell=2k$ is even.
After we multiply out Fourier expansions and extract the constant term as above, the constant term takes the shape
\[
\sum_{r \geq 1} \frac{q^{2r} r^s p_{2k}(q^r)}{(1-q^r)^{2k+1}} - (-1)^{s} \sum_{r \geq 1} \frac{q^{r} r^s p_{2k}(q^r)}{(1-q^r)^{2k+1}},
\]
where $s:=m - 2 n-1$ could be either even or odd. For $s$ even, we have
\begin{equation*} 
\sum_{r \geq 1} \frac{r^s q^{2r} p_{2k}(q^r)}{(1-q^r)^{2k+1}} - (-1)^{s} \sum_{r \geq 1} \frac{r^s q^{r} p_{2k}(q^r)}{(1-q^r)^{2k+1}}=- \sum_{r \geq 1} \frac{r^s q^{r}p_{2k}(q^r)}{(1-q^r)^{2k}}
\end{equation*}
whereas for $s$ odd, we obtain
\begin{equation*} 
\sum_{r \geq 1} \frac{r^s q^{2r} p_{2k} (q^r)}{(1-q^r)^{2k+1}} + \sum_{r \geq 1} \frac{r^s q^{r} p_{2k}(q^r)}{(1-q^r)^{2k+1}}=\sum_{r \geq 1} \frac{r^s q^{r}(1+q^r) p_{2k}(q^r)}{(1-q^r)^{2k+1}}.
\end{equation*}
By Lemma \ref{euler}, both series are in $\mathcal{QM}$.
Similarly we argue for $\ell$ odd. That completes the proof.
\end{proof}

Denote by $D_z=\frac{1}{2 \pi {\rm i} } \frac{\partial}{\partial z}$ and ${D}_\tau=\frac{1}{2 \pi {\rm i}} \frac{\partial}{\partial \tau},$ and let $s(x,y) \in \mathbb{C}[x,y]$ be a polynomial with complex coefficients.
\begin{Corollary} \label{const} Let $f(z;\tau) \in \mathcal{QJ}^{0}$. Then
\[
{\rm CT}_\zeta (s(D_z,D_\tau) f(z;\tau) ) \in \mathcal{QM}.
\]
\end{Corollary}
\begin{proof} We only have to argue that the ring $\mathcal{QJ}^{0}$ is closed under $D_z$ and $D_\tau$. This can be verified explicitly on the generators $J_1$, $\wp$, $\wp'$, $E_4$, $E_6$ (see, for instance, \cite{Ob, Ittersum}).
\end{proof}

An important special case involves products of derivatives of the $P$-functions introduced earlier.
A similar result appeared previously in \cite[Proposition~5.9]{GM}, though our proof differs.

\begin{Proposition} \label{P-power} We have
\begin{itemize}\itemsep=0pt
\item[$(i)$] For $m \geq 1$ and $k \geq 1$,
\[
{\rm CT}_{\zeta} \bigl(P_k(\zeta;\tau)^m \bigr) \in \mathcal{QM}.
\]
\item[$(ii)$] More generally, for every polynomial $s(x,y)$ and $k \geq 1$,
\[
{\rm CT}_{\zeta} \bigl\{ (s(D_\tau,D_{z}) P_k(\zeta;\tau))^{m} \bigr\} \in \mathcal{QM}.
\]
\end{itemize}
\end{Proposition}
\begin{proof} Part (ii) follows immediately from Theorem~\ref{quasi}, Corollary \ref{const}, and formulas (\ref{rho1})--(\ref{rhok}).
\end{proof}

\subsection[Characters of genus zero class-S vertex operator algebras]{Characters of genus zero class-$\boldsymbol{\mathcal{S}}$ vertex operator algebras}\label{section3.3}

In this section, $\mathfrak{g}$ denotes a finite dimensional Lie algebra of ADE type for simplicity. In the introduction, we
constructed a $q$-MZV function associated with $\mathfrak{g}$, $\zeta_{\mathfrak{g},q}(k_{\alpha}; \alpha \in \Delta_+)$. Letting $k_{\alpha} \geq 1$, we consider a slightly different $q$-MZV model where all $k_{\alpha}$ are equal and
there is an extra polynomial dependence on the labels of $\lambda \in P_+$:
\[
\zeta_{\mathfrak{g},q}^{s}(k):=\sum_{\lambda \in P_+} (\dim (V(\lambda))^{s} \frac{q^{k \langle\lambda+\rho,\rho\rangle}}{\prod_{\alpha \in \Delta_+} (1-q^{\langle \lambda, \alpha+\rho\rangle})^{k}}.
\]
The boldface symbol $\mathbf{V}_{\mathfrak g,k}$ denotes Arakawa's vertex operator algebra attached to the Lie algebra~$\mathfrak g$ (actually our notation $\mathbf V_{\mathfrak g,k}$ corresponds to Arakawa's $\mathbf V^{\mathcal S}_{G,b}$ with $\mathfrak g=
\operatorname{Lie}(G)$ and $b=k+2$; see~\cite{Arakawa})
and the positive integer parameter~$k$; it should not be confused with the ordinary italic letter $k$ used as a summation or weight parameter elsewhere.
Here $V(\lambda)$ is the finite-dimensional $\mathfrak{g}$-module of highest weight $\lambda$, so that $\dim (V(\lambda))$ can be computed using the Weyl dimension formula and $\rho=\frac12 \sum_{\alpha \in \Delta_+} \alpha$ is the half-sum of positive roots. Then we have a remarkable result of Arakawa on class-$\mathcal{S}$ vertex operator algebras motivated by the previous very important work of Beem, Rastelli and others~\cite{BR,Beem}.

\begin{Theorem}[\cite{Arakawa}] For every $k \geq 3$, there is a vertex operator algebra ${\bf V}_{\mathfrak{g},k}$ such that
\begin{equation} \label{Ar-char}
{\rm ch}[{\bf V}_{\mathfrak{g},k}](\tau)=\frac{1}{\eta(\tau)^{m_{\mathfrak{g},k}}} \zeta_{\mathfrak{g},q}^{k}(k-2),
\end{equation}
where $m_{\mathfrak{g},k}=\dim (\mathfrak{g})k-(k-2){\rm rank}(\mathfrak{g})$ is a positive integer and $\eta(\tau)=q^{1/24}\prod_{i \geq 1}(1-q^i)$ is the Dedekind $\eta$ function.
\end{Theorem}

\begin{Remark}
Although this seems to be implicit in the literature, it is expected that
\[
\dim X_{\mathbf V_{\mathfrak g,k}}=m_{\mathfrak g,k},
\]
where $X_{\mathbf V_{\mathfrak g,k}}$ is the associated variety of $\mathbf V_{\mathfrak g,k}$. In the 4d/2d correspondence, the associated variety is expected to recover the Higgs branch of the corresponding four-dimensional theory, while the vertex operator algebra character is identified with its Schur index. Thus formula~(\ref{Ar-char}) suggests that the generalized $q$-MZV is the infinite-dimensional analogue of the numerator in a Hilbert-series presentation: it contains the non-free contribution after the universal product factor has been separated. We use the phrase ``Hilbert-series numerator'' here, rather than ``$h$-series'', to avoid confusion with the $h$-polynomial of a finitely generated graded algebra.
\end{Remark}

\begin{Example}[$\mathfrak{g}=\mathfrak{sl}_2$] For $k \geq 3$,
\[
\zeta_{\mathfrak{sl}_2,q}^k(k-2)=\sum_{n \geq 1} \frac{n^k q^{\frac{n}{2} (k-2)} }{(1-q^n)^{k-2}}.
\]
\end{Example}

\begin{Example}[$\mathfrak{g}=\mathfrak{sl}_3$] For $k \geq 3$,
\[
\zeta_{\mathfrak{sl}_3,q}^k(k-2)=\sum_{m_1,m_2 \geq 1} \left( m_1 m_2 (m_1+m_2) \right)^k \frac{q^{(m_1 +m_2)(k-2)}}{(1-q^{m_1})^{k-2} (1-q^{m_2})^{k-2}(1-q^{m_1+m_2})^{k-2}}.
\]
\end{Example}

It is not hard to see that we also have
\begin{align*} 
 \zeta_{\mathfrak{sl}_{n+1},q}^k(k-2)
 =\sum_{m_1,\dots,m_n \geq 1} \left( \prod_{1 \leq i \leq j \leq n} (m_i+\cdots + m_j)^k\right) \cdot \frac{q^{(k-2)( \frac12 \sum_{i=1}^n (i-1+1)(n-i+1)m_i)}}{\prod_{1 \leq i \leq j \leq n} (1-q^{m_i+\cdots + m_j})^{k-2}}.
\end{align*}
This can be rewritten as
\begin{align*} 
 \zeta_{\mathfrak{sl}_{n+1},q}^k(k-2) =\sum_{m_1,\dots,m_n \geq 1} \prod_{1 \leq i \leq j \leq n} \left( (m_i+\cdots + m_j)^k \frac{q^{\frac12(k-2)( i m_i+\cdots + j m_j)}}{(1-q^{m_i+\cdots + m_j})^{k-2}}\right).
\end{align*}

Based on the $\mathfrak{sl}_2$ case studied by Beem and Rastelli \cite{Beem}, it is expected the following.

\begin{Conjecture} \label{MLDE} We have
\begin{itemize}\itemsep=0pt
\item[$(1)$] For $k \geq 3$ even, $\zeta_{\mathfrak{sl}_{n+1},q}^k(k-2)$ belongs to $\mathcal{QM}$ and for $k \geq 3$ odd it belongs to $\mathcal{QM}[2]:=\mathbb{Q}\bigl[E_2(\tau), \theta_2(\tau)^4,\theta_3(\tau)^4\bigr]$.\footnote{It is possible to formulate a slightly stronger conjecture depending on whether or not $\rho \in Q$, where $Q$ is the root lattice of~$\mathfrak{g}$.}
\item[$(2)$]For $k \geq 3$, the character ${\rm ch}[\mathbf{V}_{\mathfrak{g},k}](\tau)$ is a solution of a modular linear differential equation $($MLDE$)$ with coefficients in $\mathcal{QM}$ or $\mathcal{QM}[2]$, depending on the parity of~$k$.
\end{itemize}
\end{Conjecture}
Beem and Rastelli verified these conjectures for numerous examples for $\mathfrak{g}=\mathfrak{sl}_2$ (for both $k$ even and odd) and also for additional examples coming from
genus $g \geq 1$ theories. These higher-genus theories are not discussed in Arakawa's paper but presumably correspond to characters of (mostly unknown) vertex operator algebras due to 2d/4d correspondences \cite{BR}. In \cite{BVM}, we proved a broader theorem that implies the even-parameter quasi-modularity assertion in Conjecture~\ref{MLDE}\,(1) for the relevant $A$-type series. It does not by itself prove the MLDE assertion in Conjecture~\ref{MLDE}\,(2). Theorem~\ref{sl3-qm} below supplies a direct, self-contained constant-term proof for the first non-trivial case $\mathfrak{g}=\mathfrak{sl}_3$, while Proposition~\ref{sl3-qm2} treats the corresponding symmetrized family.

\subsection[Bi-brackets of multiple q-zeta values of a Lie algebra g]{Bi-brackets of multiple $\boldsymbol{q}$-zeta values of a Lie algebra $\boldsymbol{\mathfrak{g}}$}\label{section3.4}

The notation in this subsection is independent of the one-variable zeta values in (\ref{qmzvg}); here the square-bracket notation records two arrays of parameters, following the terminology of bi-brackets.
Motivated by Arakawa's character formula, we can also define more general multiple $q$-zeta series, reminiscent of the {\em bi-brackets} construction in the setup of $q$-MZVs \cite{BK2,BK1,BK3}. For simplicity, we again only discuss the $A$-type here so the positive roots are parametrized by the segments $[i,j]$, $i \leq j$.

Let $k_{[i,j]} \geq 1$, $1 \leq i \leq j$ and $s_{[i,j]} \geq 0$ be integers, and
\begin{align*}
 & \zeta_{\mathfrak{sl}_{n+1},q} \left[ \begin{matrix} k_{[i,j]}; & 1 \leq i \leq j \leq n \\ s_{[i,j]}; & 1 \leq i \leq j \leq n \end{matrix} \right]\nonumber\\
 &\qquad :=\sum_{m_1,\dots,m_n \geq 1} \ \prod_{1 \leq i \leq j \leq n} (m_i+\cdots + m_j)^{s_{[i,j]}} \frac{q^{\frac12 k_{[i,j]} (i m_i+\cdots + j m_j)}}{(1-q^{m_i+\cdots + m_j})^{k_{[i,j]}}}. \label{type-An-alt}
\end{align*}
With this definition, we have
\[
\zeta_{\mathfrak{sl}_{n+1},q} \left[ \begin{matrix} k_{[i,j]} =k; & 1 \leq i \leq j \leq n \\ s_{[i,j]} =k+2; & 1 \leq i \leq j \leq n \end{matrix} \right]
=\zeta_{\mathfrak{sl}_{n+1},q}^{k+2} (k).
\]
Based on many examples, we expect a symmetric formula to hold.

\begin{Conjecture} Assume that all $k_{[i,j]}$ are even and $s=s_{[i,j]}$ for all $[i,j]$. Then
\[
\sum_{\sigma \in W} \zeta_{\mathfrak{sl}_{n+1},q} \left[ k_{\sigma[1,1]},k_{\sigma[1,2]}, \dots, k_{\sigma[n,n]};s \right] \in \mathcal{QM},
\]
where $k_{-[i,j]}=k_{[i,j]}$.
\end{Conjecture}

One can also formulate a similar conjecture for other root systems.
We also expect quasi-modularity to hold for
$\zeta_{\mathfrak{sl}_{n+1},q} \left[ \begin{smallmatrix} k_{[i,j]}; & 1 \leq i \leq j \leq n \\ s_{[i,j]}; & 1 \leq i \leq j \leq n \end{smallmatrix} \right]$
provided that $k_{\alpha}$ and $s_{\alpha}$ have the same parity.

\subsection[Example: Quasi-modularity for sl(3)]{Example: Quasi-modularity for $\boldsymbol{\mathfrak{sl}(3)}$}\label{section3.5}

The quasi-modularity of $\mathfrak{sl}(2)$ $q$-MZV values and their modifications is well known to experts, and therefore we only briefly comment on this case.
Recall the Lie-theoretic $q$-multiple zeta values introduced in
\eqref{qmzvg}. In the simplest case $\mathfrak g=\mathfrak{sl}_2$, one has
\[
\zeta_{\mathfrak{sl}_2,q}(2k)
=
\sum_{n\ge1}\frac{q^{nk}}{(1-q^n)^{2k}}.
\]
It is well known that, up to an additive constant, these Lambert series are
linear combinations of Eisenstein series and consequently,
\[
\zeta_{\mathfrak{sl}_2,q}(2k)\in\mathcal{QM}.
\]
Similarly, for every $k\geq1$, there exist constants
$c_{k,i}\in\mathbb{Q}$, $1\leq i\leq k$, such that, for every
$n\geq1$ and $a\geq0$, the following identity of rational functions in $q$
holds:
\begin{equation}\label{lin-comb}
\frac{n^{2k+a}q^{nk}}{(1-q^n)^{2k}}
=
\sum_{i=1}^{k}
c_{k,i}
\left(q\frac{{\rm d}}{{\rm d}q}\right)^{2k+1-2i}
\frac{n^{2i-1+a}q^n}{1-q^n}.
\end{equation}

Summing over $n$ and letting $a=2$ immediately gives
\[
\zeta_{\mathfrak{sl}_2,q}^{2k+2}(2k)=\sum_{i=1}^{k} c_{k,i} \left(q \frac{{\rm d}}{{\rm d}q}\right)^{2k+1-2i} E_{2+2i}(q) \in \mathbb{Q}[E_2,E_4,E_6].
\]

\begin{Corollary}
For integers $k\geq1$ and $a\geq0$, we have
\begin{align}
&\sum_{i=1}^{k} c_{k,i}
D_\tau^{\,2k+1-2i}
\left(
\sum_{n\geq1}
\left(
\frac{n^{2i-1+a}q^n\zeta^n}{1-q^n}
+
(-1)^{a+1}
\frac{n^{2i-1+a}q^n\zeta^{-n}}{1-q^n}
\right)
\right)\nonumber
\\
&\qquad =
\sum_{n\geq1}
\left(
\frac{n^{2k+a}q^{nk}\zeta^n}{(1-q^n)^{2k}}
+
(-1)^a
\frac{n^{2k+a}q^{nk}\zeta^{-n}}{(1-q^n)^{2k}}
\right).\label{form-sl3}
\end{align}
\end{Corollary}
Equipped with (\ref{form-sl3}) and results from Section~\ref{section3.3}, we are now ready to prove the following.

\begin{Theorem} \label{sl3-qm} Let $k\geq1$ and let $a\geq0$ be even. Then
\[
\sum_{m_1,m_2\geq1}\left(m_1m_2(m_1+m_2)\right)^{2k+a}
\frac{q^{2k(m_1+m_2)}}{(1-q^{m_1})^{2k}(1-q^{m_2})^{2k}(1-q^{m_1+m_2})^{2k}}
\in\mathcal{QM}.
\]
In particular, taking $a=2$ gives $\zeta_{\mathfrak{sl}_3,q}^{2k+2}(2k)\in\mathcal{QM}$.
\end{Theorem}
\begin{proof} Keeping $a$ even, we prove that
\[
\zeta_{\mathfrak{sl}_3,q}^{2k+a}(2k)=\sum_{m_1,m_2 \geq 1} \left( m_1 m_2 (m_1+m_2) \right)^{2k+a} \frac{q^{2k(m_1 +m_2)}}{(1-q^{m_1})^{2k} (1-q^{m_2})^{2k}(1-q^{m_1+m_2})^{2k}}
\]
also belongs to $\mathcal{QM}$.
Keeping $a$ even, observe an identity
\begin{align*}
& {\rm CT}_{\zeta} \left(\sum_{n \geq 1} \frac{n^{2k+a} q^{nk} \zeta^{n}}{(1-q^n)^{2k}}+(-1)^a \frac{n^{2k+a} q^{n k} \zeta^{-n}}{(1-q^n)^{2k}} \right)^3 \\
&\qquad = {\rm CT}_{\zeta} \left(\sum_{n \geq 1} \frac{n^{2k+a} q^{nk} \zeta^{n}}{(1-q^n)^{2k}}+ \frac{n^{2k+a} q^{n k} \zeta^{-n}}{(1-q^n)^{2k}} \right)^3 \\
& \qquad =6 \sum_{m_1,m_2 \geq 1} \left( m_1 m_2 (m_1+m_2) \right)^{2k+a} \frac{q^{2k(m_1 +m_2)}}{(1-q^{m_1})^{2k} (1-q^{m_2})^{2k}(1-q^{m_1+m_2})^{2k}},
\end{align*}
where in the last step we expanded the triple product as triple sums and we extract the constant term from each individual term.
Formula~(\ref{form-sl3}), together with Proposition~\ref{P-power}, now gives the assertion.
\end{proof}

The value of the present argument is not merely the abstract quasi-modularity conclusion: it identifies the relevant $q$-series explicitly as constant terms of concrete index-zero quasi-Jacobi forms and gives a direct proof that can be adapted to further graph and Lie-theoretic series. This complements the more general result of \cite{BVM} and explains the mechanism behind the first non-trivial $\mathfrak{sl}_3$ examples.

The method in Theorem~\ref{sl3-qm}, together with Proposition~\ref{P-power}, can be used to prove the quasi-modularity property for additional $q$-MZVs associated with $\mathfrak{sl}_3$.
For example, expanding the triple product
\begin{align*}
& \sum_{m_1,m_2,m_3 \geq 1} \left(\frac{m_1^{2k+1} \zeta^{m_1} q^{m_1}}{1-q^{m_1}} +\frac{m_1^{2k+1} \zeta^{-m_1}}{1-q^{m_1}} \right) \\
& \qquad {} \times \left(\frac{m_2^{2k+1} \zeta^{m_2} q^{m_2}}{1-q^{m_2}} +\frac{m_2^{2k+1} \zeta^{-m_2}}{1-q^{m_2}} \right)
 \left(\frac{m_3^{2k+1} \zeta^{m_3} q^{m_3}}{1-q^{m_3}} +\frac{m_3^{2k+1} \zeta^{-m_3}}{1-q^{m_3}} \right)
 \end{align*}
and extracting the constant term gives
\[
\sum_{m_1,m_2 \geq 1} \frac{\left(m_1 m_2 (m_1+m_2)\right)^{2k+1} q^{m_1+m_2} }{(1-q^{m_1})(1-q^{m_2})(1-q^{m_1+m_2})} \in \mathcal{QM}.
\]
More explicitly, for $k=1$, we get
\begin{align*}
& \sum_{m_1,m_2 \geq 1} \left( m_1 m_2 (m_1+m_2) \right)^3 \frac{q^{m_1 +m_2}}{(1-q^{m_1})(1-q^{m_2})(1-q^{m_1+m_2})} \\
& \qquad = \frac{1}{570240} E_6^2(\tau)+\frac{1}{798336} E_4(\tau)^3-\frac{1}{332640} E_2(\tau)E_4(\tau)E_6(\tau).
\end{align*}
Our next result concerns $\zeta_{\mathfrak{sl}_3,q}(2k)$ values.
\begin{Proposition} \label{sl3-qm2} For every $k \geq 1$,
\[
\zeta_{\mathfrak{sl}_3,q}(2k) \in \mathcal{QM}.
\]
More generally,
\begin{equation} \label{s3}
\sum_{ \sigma \in S_3} \zeta_{\mathfrak{sl}_3,q}(2k_{\sigma{(1)}},2k_{\sigma{(2)}},2k_{\sigma{(3)}}) \in \mathcal{QM}.
\end{equation}
\end{Proposition}
\begin{proof}
Recall the series $J_k(z;\tau)$ from Section~\ref{section3.1}.
Using formula (\ref{lin-comb}), up to an additive constant, we can write (here $k \geq 1$)
\[
Q_k(z;q)=\sum_{n \geq 1} \left( \frac{q^{kn}z^{n}}{(1-q^n)^{2k}}+ \frac{q^{kn}z^{-n}}{(1-q^n)^{2k}} \right)
\]
as a linear combination of $\tilde{J}_{2k}:=\frac{1}{2k}(B_{2k}-J_{2k}(z;\tau))$, which are already in $\mathcal{QJ}^0$. For instance,
\[
\sum_{n \geq 1} \left( \frac{q^{3n}z^{n}}{(1-q^n)^{6}}+ \frac{q^{3n}z^{-n}}{(1-q^n)^{6}} \right)=\frac{1}{120} \tilde{J}_6(z;\tau)-\frac{1}{24} \tilde{J}_4(z;\tau)+\frac{1}{30} \tilde{J}_2(z;\tau).
\]
The rest follows as before by expanding the triple product
\[
{\rm CT}_\zeta \left(\sum_{n \geq 1} \frac{q^{kn}\zeta^{n}}{(1-q^n)^{2k}}+ \frac{q^{kn}\zeta^{-n}}{(1-q^n)^{2k}} \right)^3=6 \zeta_{\mathfrak{sl}_3,q}(2k).
\]
To prove (\ref{s3}), we expand the triple product
\begin{align*}
& \left(\sum_{n \geq 1} \frac{q^{k_1n}\zeta^{n}}{(1-q^n)^{2k_1}}+ \frac{q^{k_1n}\zeta^{-n}}{(1-q^n)^{2k_1}} \right)\cdot \left(\sum_{n \geq 1} \frac{q^{k_2n}\zeta^{n}}{(1-q^n)^{2k_2}}+ \frac{q^{k_2n}\zeta^{-n}}{(1-q^n)^{2k_2}} \right)\\
& \qquad \times \left(\sum_{n \geq 1} \frac{q^{k_3n}\zeta^{n}}{(1-q^n)^{2k_3}}+ \frac{q^{k_3n}\zeta^{-n}}{(1-q^n)^{2k_3}} \right),
\end{align*}
which results in eight triple sums. Two sums do not contribute to the constant term, and the remaining six sums give the left-hand side of (\ref{s3}) after taking ${\rm CT}_\zeta$.
\end{proof}

\begin{Remark} We expect $\zeta_{q,\mathfrak{sl}_3}^a(k) \notin \mathcal{QM}$ whenever $a+k$ is odd. For instance, some
manipulations with $q$-series as in \cite{Andrews} gives
\begin{align*}
& \sum_{n,m \geq 1} \frac{q^{n+m}}{(1-q^n)(1-q^m)(1-q^{n+m})}
= -\left(q \frac{{\rm d}}{{\rm d}q}\right) \sum_{n \geq 1} \frac{q^n}{1-q^n}+\sum_{n \geq 1} \frac{n^2 q^n}{1-q^n}
\end{align*}
and
\[
\sum_{n,m \geq 1} \frac{n m q^{n+m}}{(1-q^n)(1-q^m)(1-q^{n+m})}=-\left(q \frac{{\rm d}}{{\rm d}q}\right)^2 \sum_{n \geq 1} \frac{q^n}{1-q^n}+\left(q \frac{{\rm d}}{{\rm d}q}\right) \sum_{n \geq 1} \frac{n^2 q^n}{1-q^n}.
\]
Sums like $\sum_{n \geq 1} \frac{q^n}{1-q^n}$ and $\sum_{n \geq 1} \frac{n^2 q^n}{1-q^n}$, including their derivatives, do not belong to $\mathcal{QM}$.
However, many series of this type are {\em holomorphic quantum modular forms} (after Zagier). Closely related $q$-series also appeared in the study of graph series of type~$A$~\cite{BJM}.
\end{Remark}

\section[Characters and supercharacters of the U family of vertex operator (super)algebras]{Characters and supercharacters of the $\boldsymbol{\mathcal{U}}$ family\\ of vertex operator (super)algebras}\label{section4}

In \cite{AM}, Adamovi\'c and the author introduced an infinite family of vertex operator (super)algebras, somewhat analogous to those appearing in the class-$\mathcal{S}$ and in the Cvitanovi\'c--Deligne series \cite{AK,BR}.
We refer to this family as the $\mathcal{U}$ family of vertex operator (super)algebras and denote its members by $\mathcal{U}^{(m)}$. As demonstrated in \cite{AM},
for each integer $m \geq 3$, $\mathcal{U}^{(m)}$ is constructed as an infinite (super)extension of the simple affine vertex operator algebra $V_{-1}(\mathfrak{sl}(m))$ of central charge $-(m+1)$.

When $m$ is odd, $\mathcal{U}^{(m)}$ is a $\mathbb{Z}$-graded vertex operator superalgebra, and it is therefore natural to consider its supercharacter ${\rm sch}[\mathcal{U}^{(m)}](\tau)$.
Recall that the \emph{supercharacter} of a vertex operator superalgebra is
defined by taking the difference of the graded characters of its even and odd
subspaces. Equivalently, it is obtained by replacing the ordinary trace in the
definition of the character with the supertrace.
When $m$ is even, $\mathcal{U}^{(m)}$ is a vertex operator algebra with $\tfrac{1}{2}\mathbb{Z}_{\geq 0}$-grading.
For simplicity, we restrict our attention here to the case where $m$ is odd; see Remark~\ref{char-ord}.
The main motivation for the study undertaken in~\cite{AM} was the following:

\begin{Conjecture}[\cite{AM}] For every $m \geq 3$, the vertex operator $($super$)$algebra $\mathcal{U}^{(m)}$ is quasi-lisse.
\end{Conjecture}

One of the main results of \cite{AM} is that the characters and supercharacters of a distinguished family of $\mathcal{U}^{(m)}$-modules (including $\mathcal{U}^{(m)}$ itself) exhibit quasi-modular behavior.
This observation provided one of the primary motivations for the conjecture.

In this section, we present two new results on supercharacters that refine and extend the results of \cite{AM}:
\begin{enumerate}\itemsep=0pt
\item We derive a more explicit formula in terms of the constant terms of powers of the Weierstrass elliptic function, in the spirit of the computations carried out in Section~\ref{section3.3}.

\item We prove that the supercharacter can be expressed in terms of quasi-modular forms of a~specific type on the congruence subgroup $\Gamma_0(m)$.
\end{enumerate}
An important role in our construction is played by the index-zero Jacobi form
\begin{equation*} 
F_m(z;\tau)=\epsilon_m \eta(\tau)^{m+1} \frac{\vartheta(mz;m\tau)}{\vartheta(z;\tau)^m},
\end{equation*}
where $\epsilon_m$ is some non-zero constant defined in \cite{AM}. Then we have the following result.

\begin{Theorem}[\cite{AM}] For $m$ odd, we have $($in a suitable $\zeta$-range$)$
\[
{\rm CT}_{\zeta} F_m(z;\tau)={\rm sch}[\mathcal{U}^{(m)}](\tau),
\]
where as before $\zeta={\rm e}^{2 \pi {\rm i} z}$.
\end{Theorem}

\subsection[Supercharacters of U\^{}\{(m)\} via elliptic functions]{Supercharacters of $\boldsymbol{\mathcal{U}^{(m)}}$ via elliptic functions}\label{section4.1}

Using the notation from \cite{Lang}, let
\[
G_{2,N,(0,\frac{r}{N})}(\tau):=\wp \left(\frac{r}{N};\tau \right)
\]
with $1 \leq r \leq N-1$, a certain weight two Eisenstein series of level $N$.
Then we have for $A=\left(\begin{smallmatrix} a & b \\ c & d \end{smallmatrix}\right) \in {\rm SL}(2,\mathbb{Z})$,
\[
G_{2,N,(0,\frac{r}{N})}\left(\frac{a \tau+b}{c \tau +d}\right)=(c\tau+d)^2 G_{2,N,A(0,\frac{r}{N}) }(\tau).
\]
Consider now $c | N$, so that $\gcd(d,N)=1$. Then the right-hand side equals
\[
(c\tau+d)^2 G_{2,N, (\frac{cr}{N}, \frac{dr}{N})}(\tau)= (c\tau+d)^2 G_{2,N, (0, \frac{dr}{N})}(\tau).
\]
In other words, elements of $\Gamma_0(N)$ permute the $G_2$-series. Using $\wp(-z;\tau)=\wp(z;\tau)$ and letting~$N$ be
odd, we see that
\begin{equation*} 
\mathbb{G}_1(\tau):=\sum_{a=1}^{(N-1)/2} G_{2,N, (0, \frac{a}{N})}(\tau)
\end{equation*}
is invariant under $\Gamma_0(N)$. Similarly, we evaluate elementary symmetric polynomials $e_i(x_1,\dots,\allowbreak x_{(N-1)/2})$, $ 1 \leq i \leq (N-1)/2$,
at $G_{2,N,(0,\frac{r}{N})}(\tau)$ and construct genuine modular forms on $\Gamma_0(N)$, which we denote by $\mathbb{G}_i(\tau)$ (observe that
$\mathbb{G}_1(\tau)$ comes from $x_1+x_2+\cdots+x_{(N-1)/2}$).
These are modular forms of weight $2i$.

Now we are ready to give a new formula for the supercharacter of $\mathcal{U}$ expressed using
${\rm CT}_\zeta \wp(z)^k$. \begin{Proposition} \label{jacobi} We have the supercharacter formula
\[
{\rm sch}[\mathcal{U}^{(m)}](\tau)=\mu_m \frac{\eta(m \tau)^3}{\eta(\tau)^{2m-1}} \sum_{k=0}^{(m-1)/2} (-1)^k \mathbb{G}_k(\tau) {\rm CT}_\zeta \wp(z;\tau)^k,
\]
where $\mu_m \neq 0$ is a constant depending only on $m$.
\end{Proposition}
\begin{proof}
We first recall identity (\ref{rho-sigma}) that allows us to write
\begin{align*}
& \wp(z;\tau)-\wp\left(\frac{j}{m};\tau\right)= \frac{\theta\bigl(z+\frac{j}{m};\tau\bigr) \theta(z-\frac{j}{m};\tau)}{ 4 \pi^2 \theta(z;\tau)^2 \theta\bigl(\frac{j}{m};\tau\bigr)^2 \eta(\tau)^6}.
\end{align*}
Multiplying out these relations from $j=1,\dots,\frac{m-1}{2}$ (here $m=2k+1$ is odd!), using
\[
\bigl(1-x^{2k+1}\bigr)=(1-x)\prod_{j=1}^{k}\bigl(1-x {\rm e}^{\frac{2 \pi j }{2k+1}}\bigr)\prod_{j=1}^{k}\bigl(1-x {\rm e}^{-\frac{2 \pi j }{2k+1}}\bigr),
\]
 modulo Euler factors, we can express the Jacobi form $F_m(z;\tau)$ as the product of $\wp(\tau)-\wp\bigl(\frac{k}{m}; \tau\bigr)$, where $1 \leq k \leq (m-1)/2$. After collecting the $z$-independent theta and eta factors, put $r=(m-1)/2$ and $a_j=\wp(j/m;\tau)$. Vieta's formula gives
\[
\prod_{j=1}^{r}(X-a_j)=\sum_{k=0}^{r}(-1)^{r-k}e_{r-k}(a_1,\dots,a_r)X^k.
\]
Taking $X=\wp(z;\tau)$ shows explicitly that the coefficient of $\wp(z;\tau)^k$ is $(-1)^{r-k}e_{r-k}(a_1,\dots,a_r)$; these elementary symmetric functions are precisely the forms denoted by $\mathbb{G}_{r-k}(\tau)$, up to the common normalization absorbed into~$\mu_m$.
\end{proof}

As an application, we compute the supercharacter for $m=3$.
We only have to examine $\wp(z;\tau)-\wp(\frac{1}{3},\tau)$.
Explicit computation with Eisenstein series gives
\[
\wp\left(\frac13;\tau\right)=\frac32G_2(\tau)-\frac{9}{2}G_2(3 \tau).
\]
Combined with the fact that the constant term of $\wp(z;\tau)$ is $G_2(\tau)$, we obtain
\[
G_2(\tau)-(\frac32G_2(\tau)-\frac{9}{2}G_2(3 \tau) )=-\frac12 G_2(\tau)+\frac92 G_2(3 \tau).
\]
This agrees with the formula obtained in \cite{AM}:
\begin{align*}
{\rm sch}[\mathcal{U}^{(3)}](\tau)&=\frac{\eta^{3}(3 \tau)}{\eta^5(\tau)}\left(-\frac{1}{8}E_{2}(\tau)+\frac98 E_{2}(3 \tau)\right) \\
&=q^{\frac16}\bigl(1+8q+44q^2+152q^3+487q^4+1352q^5+\cdots\bigr).
\end{align*}
We should point out that $E_{2}(\tau)-9 E_{2}(3 \tau)$ is not modular on $\Gamma_0(3)$ but quasi-modular. More precisely, up to a scalar, it equals the logarithmic
derivative of $\frac{\eta(3 \tau)^3}{\eta(\tau)}$.
Another slightly more involved example is $m=5$. Analyzing $\wp\bigl(\frac{i}{5};\tau\bigr)$, $ 1\leq i \leq 2$, gives
\begin{align*}
&{\rm sch}[\mathcal{U}^{(5)}](\tau) \\
&\quad= \frac{\eta^{3}(5 \tau)}{\eta^9(\tau)}\left( \frac{25}{1152} E^2_{2}(\tau)- \frac{125}{192} E_{2}(\tau) E_2(5 \tau) +\frac{3125}{1152} E_2(5 \tau)^2 +\frac{1}{576}
E_4(\tau)-\frac{625}{576}E_4(5 \tau) \right) \\
 & \quad =q^{\frac23}\bigl(1 + 24 q + 249 q^2 + 1750 q^3 + 9750 q^4 + 45750 q^5+ \cdots \bigr).
\end{align*}
The appearance of the series $E_{2k}(m \tau)$ in these formulas is not an accident.
\begin{Theorem} \label{charum} Assume $m\geq3$ is odd. Then we have
\[
{\rm sch}[\mathcal{U}^{(m)}](\tau)= \frac{\eta(m \tau)^3}{\eta(\tau)^{2m-1}} \mathbb{F}_m(\tau),
\]
where $\mathbb{F}_m(\tau) \in \mathbb{Q}[E_2(\tau),E_4(\tau),E_6(\tau),E_2(m \tau),E_4(m \tau),E_6(m \tau)]$ is a quasi-modular form on $\Gamma_0(m)$ of weight $m-1$.
\end{Theorem}
\begin{proof} Let ${\mathcal{QM}}(m):=\mathbb{Q}[E_2(\tau),E_4(\tau),E_6(\tau),E_2(m \tau),E_4(m \tau),E_6(m \tau)]$.
We only have to find a Laurent expansion that involves suitable Eisenstein series.
For this, we use a well-known identity \cite{Lang}
\[
\vartheta(z;\tau)=-2 \pi z \eta(\tau)^3 \cdot \exp \left(-\sum_{k=1}^\infty \frac{G_{2k}(\tau)}{2k} z^{2k} \right).
\]
The Laurent coefficients in the expansion around $z=0$ of
\[
\frac{\vartheta(mz; m\tau)}{\vartheta(z; \tau)^m}
\]
are of the form $E_{2i}(\tau)$ and $E_{2i}(m\tau)$. Consequently, they lie in the ring generated by $E_{2i}(\tau)$ and $E_{2i}(m\tau)$ for $i=1,2,3$, and hence belong to $\mathcal{QM}(m)$.

Using the same argument as in Theorem~\ref{quasi} or in \cite{AM}, we then show that the constant term can be expressed in terms of these Eisenstein series together with ordinary quasi-modular forms.

Finally, to see that $\mathbb{F}_m(\tau)$ is a quasi-modular form on the congruence subgroup $\Gamma_0(m)$, it suffices to observe that
\[
E^{(m)}_2(\tau) := E_2(\tau) - m E_2(m\tau), \qquad
E_4(m\tau), \qquad
E_6(m\tau)
\]
are modular forms on $\Gamma_0(m)$.
Therefore, any homogeneous polynomial in
\[
E_2(\tau),\ E_4(\tau),\ E_6(\tau),\ E_2(m\tau),\ E_4(m\tau),\ \text{and} \ E_6(m\tau)
\]
is a quasi-modular form on $\Gamma_0(m)$.
\end{proof}

\begin{Remark} \label{char-ord} The same method should apply to the ordinary characters of
$\mathcal{U}^{(m)}$ and its irreducible modules for both odd and even $m$.
The present paper treats only the supercharacters in the odd case, where they
have a particularly nice form. More precisely, using the formulas from~\cite{AM}, we have
\[
\operatorname{ch}[\mathcal{U}^{(m)}](\tau) = \epsilon_m \eta(\tau)^{m+1} \operatorname{CT}_\zeta \left( \frac{\vartheta\left(zm + \frac{1}{2}; m\tau\right)}{\vartheta(z; \tau)^m} \right).
\]
for some non-zero number $\epsilon_m$. Hence, as in Theorem~\ref{charum} we obtain the same type of Eisenstein series -- however the formula is now more difficult to write down due to the shift in the $z$-variable.
\end{Remark}

\section{Final remarks}\label{section5}

We end with a few comments.
\begin{itemize}\itemsep=0pt
\item[(i)] We are not aware of examples of graphs $\Gamma$ for which $q^a H_{\Gamma}(q)$ is modular except if $\Gamma$ is totally disconnected (no edges). We are also not aware of non-isomorphic {\em connected} graphs whose graph series are identical.
However, of course, there are edge ideals of connected graphs with the same Hilbert series. For instance, $R_{C_5}$ and $R_{\Gamma}$, where $I(\Gamma)=\langle x_1 x_2,x_1 x_3,x_3 x_4,x_4 x_5,x_5 x_3 \rangle$ have the same Hilbert series, although their graphs are non-isomorphic.

To the best of our knowledge, the graphs in Section~\ref{section2} are the only infinite families of connected graphs (up to isomorphism) whose graph series can be reduced to familiar $q$-series. Even elementary graphs such as paths (and especially cycles) are quite non-trivial to analyze \cite{BJM}.

\item[(ii)] Based on the $\mathfrak{sl}(3)$ example, we believe that ${\rm ch}[\mathbf{V}_{\mathfrak{g},k}]$ appears as the constant term in the Fourier expansion (in appropriate range) of a nice family of meromorphic quasi-Jacobi forms of index zero in several variables. At the same time, there is computational evidence that characters ${\rm ch}[\mathbf{V}_{\mathfrak{g},k}]$ can be interpreted as certain
$1$-point correlation functions on the torus:
\[
{\rm tr}_{\mathcal{F}} \ o(a) q^{L(0)-\frac{c}{24}}
\]
computed on a suitable bosonic space $\mathcal{F}$ or a $\mathbb{Z}_2$-twisted $\mathcal{F}$-module, for particular choice of $a \in \mathcal{F}$. If so, then \cite{DMN} would immediately imply Conjecture~\ref{MLDE}\,(i).

\item[(iii)] The approach based on $1$-point functions, combined with \cite{M2, M1}, can be used to give another proof of quasi-modularity of ``correlation functions'' in the theory of Hilbert schemes of point on surfaces \cite{Carlsson} and various Bloch--Okounkov type formulas \cite{Ittersum}.
Another interesting aspect of $q$-MZVs is in the theory of Hilbert schemes of points; see~\cite{Okounkov} and especially~\cite{Qin}. It is unclear to us whether their work is connected with ours.
\end{itemize}

\subsection*{Acknowledgements}

The author received partial support from the Collaboration Grant for Mathematicians 709563 awarded by the Simons Foundation and the NSF grant DMS-2101844.
This work is funded by the EU NextGeneration under the Juraj Dobrila University
of Pula institutional research project number IIP\_{}UNIPU\_{}010159.

 Results from this paper were presented at various seminars and also at the Banff meeting on ``Quantum Field Theories and Quantum Topology Beyond Semisimplicity'',
 October--November 2021. We would like to express our gratitude to the organizers for the invitation. After this paper was posted on the arXiv, we were informed by Palash Singh about two very recent papers
\cite{Beem-new, Pan} dealing with indices of class-$\mathcal{S}$ theories.

We are grateful to the anonymous referees for their careful reading and numerous constructive suggestions, which substantially improved the exposition of this paper.


\pdfbookmark[1]{References}{ref}
\LastPageEnding


\begin{thebibliography}{99}
\footnotesize\itemsep=0pt

\bibitem{AM}
Adamovi\'{c} D., Milas A., On some vertex algebras related to
 {$V_{-1}(\mathfrak{sl}(n))$} and their characters,
 \href{https://doi.org/10.1007/s00031-020-09617-w}{\textit{Transform. Groups}}
 \textbf{26} (2021), 1--30,
 \href{https://arxiv.org/abs/1805.09771}{arXiv:1805.09771}.

\bibitem{Andrews0}
Andrews G.E., The theory of partitions, \textit{Encyclopedia Math. Appl.},
 Vol.~2, Addison-Wesley Publishing Co., Reading, Mass., 1976.

\bibitem{Andrews}
Andrews G.E., Stacked lattice boxes,
 \href{https://doi.org/10.1007/BF01608779}{\textit{Ann. Comb.}} \textbf{3}
 (1999), 115--130.

\bibitem{AR}
Andrews G.E., Rose S.C.F., Mac{M}ahon's sum-of-divisors functions, {C}hebyshev
 polynomials, and quasi-modular forms,
 \href{https://doi.org/10.1515/crelle.2011.179}{\textit{J.~Reine Angew.
 Math.}} \textbf{676} (2013), 97--103,
 \href{https://arxiv.org/abs/1010.5769}{arXiv:1010.5769}.

\bibitem{Arakawa}
Arakawa T., Chiral algebras of class $\mathcal{S}$ and {M}oore--{T}achikawa
 symplectic varieties,
 \href{https://arxiv.org/abs/1811.01577}{arXiv:1811.01577}.

\bibitem{AK}
Arakawa T., Kawasetsu K., Quasi-lisse vertex algebras and modular linear
 differential equations, in Lie Groups, Geometry, and Representation Theory,
 \textit{Progr. Math.}, Vol. 326,
 \href{https://doi.org/10.1007/978-3-030-02191-7_2}{Birkh\"auser}, Cham, 2018,
 41--57, \href{https://arxiv.org/abs/1610.05865}{arXiv:1610.05865}.

\bibitem{BK2}
Bachmann H., K\"uhn U., A short note on a conjecture of {O}kounkov about a
 $q$-analogue of multiple zeta values,
 \href{https://arxiv.org/abs/1407.6796}{arXiv:1407.6796}.

\bibitem{BK1}
Bachmann H., K\"uhn U., The algebra of generating functions for multiple
 divisor sums and applications to multiple zeta values,
 \href{https://doi.org/10.1007/s11139-015-9707-7}{\textit{Ramanujan~J.}}
 \textbf{40} (2016), 605--648,
 \href{https://arxiv.org/abs/1309.3920}{arXiv:1309.3920}.

\bibitem{BK3}
Bachmann H., K\"uhn U., A dimension conjecture for $q$-analogues of multiple
 zeta values, in Periods in Quantum Field Theory and Arithmetic,
 \textit{Springer Proc. Math. Stat.}, Vol.~314,
 \href{https://doi.org/10.1007/978-3-030-37031-2_9}{Springer}, Cham, 2020,
 237--258, \href{https://arxiv.org/abs/1708.07464}{arXiv:1708.07464}.

\bibitem{BR}
Beem C., Lemos M., Liendo P., Peelaers W., Rastelli L., van Rees B.C., Infinite
 chiral symmetry in four dimensions,
 \href{https://doi.org/10.1007/s00220-014-2272-x}{\textit{Comm. Math. Phys.}}
 \textbf{336} (2015), 1359--1433,
 \href{https://arxiv.org/abs/1312.5344}{arXiv:1312.5344}.

\bibitem{Beem}
Beem C., Rastelli L., Vertex operator algebras, {H}iggs branches, and modular
 differential equations,
 \href{https://doi.org/10.1007/jhep08(2018)114}{\textit{J.~High Energy Phys.}}
 \textbf{2018} (2018), no.~8, 114, 72~pages,
 \href{https://arxiv.org/abs/1707.07679}{arXiv:1707.07679}.

\bibitem{Beem-new}
Beem C., Singh P., Razamat S.S., Schur indices of class {$\mathcal{S}$} and
 quasimodular forms,
 \href{https://doi.org/10.1103/physrevd.105.085009}{\textit{Phys. Rev.~D}}
 \textbf{105} (2022), 085009, 15~pages,
 \href{https://arxiv.org/abs/2112.10715}{arXiv:2112.10715}.

\bibitem{Bradley}
Bradley D.M., Multiple {$q$}-zeta values,
 \href{https://doi.org/10.1016/j.jalgebra.2004.09.017}{\textit{J.~Algebra}}
 \textbf{283} (2005), 752--798,
 \href{https://arxiv.org/abs/math/0402093}{arXiv:math/0402093}.

\bibitem{BFM}
Bringmann K., Folsom A., Mahlburg K., Quasimodular forms and
 {$s\ell(m|m)^\wedge$} characters,
 \href{https://doi.org/10.1007/s11139-014-9621-4}{\textit{Ramanujan~J.}}
 \textbf{36} (2015), 103--116,
 \href{https://arxiv.org/abs/1302.4040}{arXiv:1302.4040}.

\bibitem{BJM}
Bringmann K., Jennings-Shaffer C., Milas A., Graph schemes, graph series, and
 modularity,
 \href{https://doi.org/10.1016/j.jcta.2023.105749}{\textit{J.~Combin. Theory
 Ser.~A}} \textbf{197} (2023), 105749, 37~pages,
 \href{https://arxiv.org/abs/2105.05660}{arXiv:2105.05660}.

\bibitem{BVM}
Bringmann K., Milas A., van Ittersum J.W.M., Quasimodularity  of $q$-series associated to root lattices, {i}n preparation.

\bibitem{CLM}
Calinescu C., Lepowsky J., Milas A., Vertex-algebraic structure of the
 principal subspaces of level one modules for the untwisted affine {L}ie
 algebras of types {$A$, $D$, $E$},
 \href{https://doi.org/10.1016/j.jalgebra.2009.09.029}{\textit{J.~Algebra}}
 \textbf{323} (2010), 167--192,
 \href{https://arxiv.org/abs/0908.4054}{arXiv:0908.4054}.

\bibitem{Carlsson}
Carlsson E., Vertex operators, {G}rassmannians, and {H}ilbert schemes,
 \href{https://doi.org/10.1007/s00220-010-1123-7}{\textit{Comm. Math. Phys.}}
 \textbf{300} (2010), 599--613,
 \href{https://arxiv.org/abs/0910.5528}{arXiv:0910.5528}.

\bibitem{DMN}
Dong C., Mason G., Nagatomo K., Quasi-modular forms and trace functions
 associated to free boson and lattice vertex operator algebras,
 \href{https://doi.org/10.1155/S1073792801000204}{\textit{Int. Math. Res.
 Not.}} \textbf{2001} (2001), 409--427,
 \href{https://arxiv.org/abs/math/0004132}{arXiv:math/0004132}.

\bibitem{DFR}
Dotsenko V., Feigin E., Reineke M., Koszul algebras and {D}onaldson--{T}homas
 invariants, \href{https://doi.org/10.1007/s11005-022-01604-4}{\textit{Lett.
 Math. Phys.}} \textbf{112} (2022), 106, 39~pages,
 \href{https://arxiv.org/abs/2111.07588}{arXiv:2111.07588}.

\bibitem{DM}
Dotsenko V., Mozgovoy S., D{T} invariants from vertex algebras,
 \href{https://doi.org/10.1017/S1474748024000288}{\textit{J.~Inst. Math.
 Jussieu}} \textbf{24} (2025), 291--339,
 \href{https://arxiv.org/abs/2108.10338}{arXiv:2108.10338}.

\bibitem{Efimov}
Efimov A.I., Cohomological {H}all algebra of a symmetric quiver,
 \href{https://doi.org/10.1112/S0010437X12000152}{\textit{Compos. Math.}}
 \textbf{148} (2012), 1133--1146,
 \href{https://arxiv.org/abs/1103.2736}{arXiv:1103.2736}.

\bibitem{Ekholm}
Ekholm T., Gruen A., Gukov S., Kucharski P., Park S., Sto\v{s}i\'{c} M.,
 Su{\l}kowski P., Branches, quivers, and ideals for knot complements,
 \href{https://doi.org/10.1016/j.geomphys.2022.104520}{\textit{J.~Geom.
 Phys.}} \textbf{177} (2022), 104520, 75~pages,
 \href{https://arxiv.org/abs/2110.13768}{arXiv:2110.13768}.

\bibitem{FS}
Feigin B., Stoyanovsky A.V., Quasi-particles models for the representations of
 {L}ie algebras and geometry of flag manifold,
 \href{https://arxiv.org/abs/hep-th/9308079}{arXiv:hep-th/9308079}.

\bibitem{Goncharov}
Goncharov A.B., Multiple polylogarithms, cyclotomy and modular complexes,
 \href{https://doi.org/10.4310/MRL.1998.v5.n4.a7}{\textit{Math. Res. Lett.}}
 \textbf{5} (1998), 497--516.

\bibitem{GM}
Goujard E., M\"{o}ller M., Counting {F}eynman-like graphs: {Q}uasimodularity
 and {S}iegel--{V}eech weight,
 \href{https://doi.org/10.4171/jems/924}{\textit{J.~Eur. Math. Soc. (JEMS)}}
 \textbf{22} (2020), 365--412,
 \href{https://arxiv.org/abs/1609.01658}{arXiv:1609.01658}.

\bibitem{CM}
Herzog J., Hibi T., Distributive lattices, bipartite graphs and {A}lexander
 duality,
 \href{https://doi.org/10.1007/s10801-005-4528-1}{\textit{J.~Algebraic
 Combin.}} \textbf{22} (2005), 289--302,
 \href{https://arxiv.org/abs/math/0307235}{arXiv:math/0307235}.

\bibitem{Hoffman}
Hoffman M.E., Algebraic aspects of multiple zeta values, in Zeta Functions,
 Topology and Quantum Physics, \textit{Dev. Math.}, Vol.~14,
 \href{https://doi.org/10.1007/0-387-24981-8_4}{Springer}, New York, 2005,
 51--73, \href{https://arxiv.org/abs/math/0309425}{arXiv:math/0309425}.

\bibitem{Zagier1}
Ihara K., Kaneko M., Zagier D., Derivation and double shuffle relations for
 multiple zeta values,
 \href{https://doi.org/10.1112/S0010437X0500182X}{\textit{Compos. Math.}}
 \textbf{142} (2006), 307--338.

\bibitem{JM1}
Jennings-Shaffer C., Milas A., On {$q$}-series identities for false theta
 series, \href{https://doi.org/10.1016/j.aim.2020.107411}{\textit{Adv. Math.}}
 \textbf{375} (2020), 107411, 22~pages,
 \href{https://arxiv.org/abs/2001.11368}{arXiv:2001.11368}.

\bibitem{JM2}
Jennings-Shaffer C., Milas A., Further {$q$}-series identities and conjectures
 relating false theta functions and characters, in Lie groups, Number Theory,
 and Vertex Algebras, \textit{Contemp. Math.}, Vol.~768,
 \href{https://doi.org/10.1090/conm/768/15467}{American Mathematical Society},
 Providence, RI, 2021, 253--269,
 \href{https://arxiv.org/abs/2005.13620}{arXiv:2005.13620}.

\bibitem{Kawasetsu}
Kawasetsu K., The free generalized vertex algebras and generalized principal
 subspaces,
 \href{https://doi.org/10.1016/j.jalgebra.2015.07.014}{\textit{J.~Algebra}}
 \textbf{444} (2015), 20--51,
 \href{https://arxiv.org/abs/1502.05276}{arXiv:1502.05276}.

\bibitem{Quivers}
Kucharski P., Reineke M., Sto\v{s}i\'{c} M., Su{\l}kowski P., Knots-quivers
 correspondence,
 \href{https://doi.org/10.4310/atmp.2019.v23.n7.a4}{\textit{Adv. Theor. Math.
 Phys.}} \textbf{23} (2019), 1849--1902,
 \href{https://arxiv.org/abs/1707.04017}{arXiv:1707.04017}.

\bibitem{Lang}
Lang S., Elliptic functions, 2nd ed., \textit{Grad. Texts in Math.}, Vol.~112,
 \href{https://doi.org/10.1007/978-1-4612-4752-4}{Springer}, New York, 1987.

\bibitem{Hao}
Li H., Some remarks on associated varieties of vertex operator superalgebras,
 \href{https://doi.org/10.1007/s40879-021-00477-6}{\textit{Eur.~J. Math.}}
 \textbf{7} (2021), 1689--1728,
 \href{https://arxiv.org/abs/2007.04522}{arXiv:2007.04522}.

\bibitem{HM}
Li H., Milas A., Jet schemes, quantum dilogarithm and {F}eigin--{S}toyanovsky's
 principal subspaces,
 \href{https://doi.org/10.1016/j.jalgebra.2023.10.033}{\textit{J.~Algebra}}
 \textbf{640} (2024), 21--58,
 \href{https://arxiv.org/abs/2010.02143}{arXiv:2010.02143}.

\bibitem{Libgober}
Libgober A., Elliptic genera, real algebraic varieties and quasi-{J}acobi
 forms, in Topology of Stratified Spaces, \textit{Math. Sci. Res. Inst.
 Publ.}, Vol.~58, Cambridge University Press, Cambridge, 2011, 95--120,
 \href{https://arxiv.org/abs/0904.1026}{arXiv:0904.1026}.

\bibitem{M2}
Milas A., Formal differential operators, vertex operator algebras and
 zeta-values,~{II},
 \href{https://doi.org/10.1016/S0022-4049(03)00036-7}{\textit{J.~Pure Appl.
 Algebra}} \textbf{183} (2003), 191--244,
 \href{https://arxiv.org/abs/math/0303154}{arXiv:math/0303154}.

\bibitem{M1}
Milas A., On certain automorphic forms associated to rational vertex operator
 algebras, in Moonshine: The First Quarter Century and Beyond, \textit{London
 Math. Soc. Lecture Note Ser.}, Vol.~372,
 \href{https://doi.org/10.1017/CBO9780511730054.016}{Cambridge University
 Press}, Cambridge, 2010, 330--357.

\bibitem{MP}
Milas A., Penn M., Lattice vertex algebras and combinatorial bases: general
 case and {$\mathcal W$}-algebras, \textit{New York J. Math.} \textbf{18}
 (2012), 621--650.

\bibitem{Ob-Serre}
Oberdieck G., A {S}erre derivative for even weight {J}acobi forms,
 \href{https://arxiv.org/abs/1209.5628}{arXiv:1209.5628}.

\bibitem{Ob}
Oberdieck G., Gromov--{W}itten invariants of the {H}ilbert schemes of points of
 a {K}3 surface, \href{https://doi.org/10.2140/gt.2018.22.323}{\textit{Geom.
 Topol.}} \textbf{22} (2018), 323--437,
 \href{https://arxiv.org/abs/1406.1139}{arXiv:1406.1139}.

\bibitem{Ohno}
Ohno Y., Okuda J., Zudilin W., Cyclic {$q$}-{MZSV} sum,
 \href{https://doi.org/10.1016/j.jnt.2011.08.001}{\textit{J.~Number Theory}}
 \textbf{132} (2012), 144--155.

\bibitem{Zagier2}
Ohno Y., Zagier D., Multiple zeta values of fixed weight, depth, and height,
 \href{https://doi.org/10.1016/S0019-3577(01)80037-9}{\textit{Indag.
 Math.~(N.S.)}} \textbf{12} (2001), 483--487.

\bibitem{Okounkov}
Okunkov A., Hilbert schemes and multiple {$q$}-zeta values,
 \href{https://doi.org/10.1007/s10688-014-0054-z}{\textit{Funct. Anal. Appl.}}
 \textbf{48} (2014), 138--144,
 \href{https://arxiv.org/abs/1404.3873}{arXiv:1404.3873}.

\bibitem{OEIS}
The on-line encyclopedia of integer sequences, \url{https://oeis.org/}.

\bibitem{Pan}
Pan Y., Peelaers W., Exact {S}chur index in closed form,
 \href{https://doi.org/10.1103/physrevd.106.045017}{\textit{Phys. Rev.~D}}
 \textbf{106} (2022), 045017, 34~pages,
 \href{https://arxiv.org/abs/2112.09705}{arXiv:2112.09705}.

\bibitem{Qin}
Qin Z., Yu F., On {O}kounkov's conjecture connecting {H}ilbert schemes of
 points and multiple {$q$}-zeta values,
 \href{https://doi.org/10.1093/imrn/rnw244}{\textit{Int. Math. Res. Not.}}
 \textbf{2018} (2018), 321--361,
 \href{https://arxiv.org/abs/1510.00837}{arXiv:1510.00837}.

\bibitem{Edge}
Renteln P., The {H}ilbert series of the face ring of a flag complex,
 \href{https://doi.org/10.1007/s003730200045}{\textit{Graphs Combin.}}
 \textbf{18} (2002), 605--619.

\bibitem{Ittersum}
van Ittersum J.W.M., Partitions and quasimodular forms: {V}ariations on the
 {B}loch--{O}kounkov theorem, Ph.D. thesis, {U}trecht {U}niversity, 2021,
 \url{https://staff.fnwi.uva.nl/j.w.m.vanittersum/files/theses/thesis_vanIttersum_withcover.pdf}.

\bibitem{Zagier}
Zagier D., Multiple zeta values, {P}reprint, 1995.

\bibitem{Zhu}
Zhu Y., Modular invariance of characters of vertex operator algebras,
 \href{https://doi.org/10.1090/S0894-0347-96-00182-8}{\textit{J.~Amer. Math.
 Soc.}} \textbf{9} (1996), 237--302.

\bibitem{Zudilin2}
Zudilin W., Algebraic relations for multiple zeta values,
 \href{https://doi.org/10.1070/RM2003v058n01ABEH000592}{\textit{Russian Math.
 Surveys}} \textbf{58} (2003), 1--29.

\bibitem{Zudilin1}
Zudilin W., Multiple $q$-zeta brackets,
 \href{https://doi.org/10.3390/math3010119}{\textit{Mathematics}} \textbf{3}
 (2015), 119--130, \href{https://arxiv.org/abs/1412.0163}{arXiv:1412.0163}.

\end{thebibliography}
\end{document}